\documentclass[12pt]{article}
\usepackage{amssymb}
\usepackage{amsmath}
\numberwithin{equation}{section}

\newcommand{\bea}{\begin{eqnarray}}
\newcommand{\eea}{\end{eqnarray}}
\newcommand{\beano}{\begin{eqnarray*}}
\newcommand{\eeano}{\end{eqnarray*}}
\newcommand{\nonu}{\nonumber \\}

\newcommand{\eps}{\epsilon}

\newcommand{\vph}{\varphi}

\newcommand{\ce}{\mbox{${\cal E}$}}
\newcommand{\cg}{\mbox{$\cal{G}$}}
\newcommand{\ch}{\mbox{$\cal{H}$}}
\newcommand{\ci}{\mbox{${\cal I}$}}
\newcommand{\cj}{\mbox{${\cal J}$}}
\newcommand{\ck}{\mbox{$\cal{K}$}}
\newcommand{\cl}{\mbox{$\cal{L}$}}

\newcommand{\ct}{\mbox{$\cal{T}$}}
\newcommand{\cu}{\mbox{${\cal U}$}}
\newcommand{\cw}{\mbox{$\cal{W}$}}

\newcommand{\wh}[1]{\widehat{#1}}

\newcommand{\mb}[1]{\ \mbox{ #1 }\ }
\newcommand{\half}{\frac{1}{2}}

\newcommand{\prf}{\underline{Proof:}\ }
\newcommand{\finprf}{\null \hfill {\rule{5pt}{5pt}}\\[2.1ex]\indent}
\newcommand{\ie}{{\it i.e.}\ }
\newtheorem{theo}{Theorem}[section]
\newtheorem{coro}[theo]{Corollary}
\newtheorem{prop}[theo]{Property}
\newtheorem{defi}[theo]{Definition}
\newtheorem{lem}[theo]{Lemma}

\newcounter{rem}
\newtheorem{rmk}[rem]{Remark}
\newcommand{\CC}{\mbox{${\mathbb C}$}}

\newcommand{\ZZ}{\mbox{${\mathbb Z}$}}
\newcommand{\NN}{\mbox{${\mathbb N}$}}

\newcommand{\II}{\mbox{${\mathbb I}$}}
\newcommand{\unity}{1 \hspace{-0.3ex} \rule{.1ex}{1.55ex}
  \rule{0.2ex}{0.08ex} \hspace{0.4ex}}

\newcommand{\FF}{\mbox{${\mathbb F}$}}
\newcommand{\JJ}{\mbox{${\mathbb J}$}}

\newcommand{\NP}[1]{Nucl.\ Phys.\ {\bf #1}}
\newcommand{\PL}[1]{Phys.\ Lett.\ {\bf #1}}

\newcommand{\CMP}[1]{Commun.\ Math.\ Phys.\ {\bf #1}}

\newcommand{\JPhys}[1]{J.\ Phys.\ {\bf #1}}
\newcommand{\IJMP}[1]{Int. J.\ Mod.\ Phys.\ {\bf #1}}

\def\opluslim{\mathop{\oplus}\limits}

\newcommand{\bara}{\bar a}
\newcommand{\barb}{\bar b}
\newcommand{\barc}{\bar c}
\newcommand{\bard}{\bar d}

\begin{document}
\renewcommand{\thefootnote}{\fnsymbol{footnote}}
\newpage
\setcounter{page}{0}

\newcommand{\LAP}{LAPTH}
\def\logo{{\bf {\huge LAPTH}}}
%
\markright{\today\dotfill SuperYangiens superW alg\dotfill DRAFT\dotfill}
\pagestyle{empty}

\centerline{\logo}

\vspace {.3cm}

\centerline{{\bf{\it\Large 
Laboratoire d'Annecy-le-Vieux de Physique Th\'eorique}}}

\centerline{\rule{12cm}{.42mm}}

\vfill

\begin{center}

  {\LARGE  {\sffamily  \cw-superalgebras as 
            truncations of super-Yangians }}\\[1cm]

\vspace{10mm}
  
{\Large C. Briot and E. Ragoucy\footnote{ragoucy@lapp.in2p3.fr}}\\[.42cm]
{\large Laboratoire de Physique Th{\'e}orique \LAP\footnote{URA 14-36 
    du CNRS, associ{\'e}e {\`a} l'Universit{\'e} de Savoie.}\\[.242cm]
    LAPP, BP 110, F-74941  Annecy-le-Vieux Cedex, France. }
\end{center}
\vfill\vfill

\begin{abstract}
        We show that some finite \cw-superalgebras based on $gl(M|N)$ 
        are truncation of the super-Yangian $Y(gl(M|N))$. In the 
        same way, we prove that finite \cw-superalgebras based on 
        $osp(M|2n)$ are truncation of the twisted super-Yangians 
        $Y(gl(M|2n))^{+}$. Using this homomorphism, we present
        these $\cw$-superalgebras in an $R$-matrix formalism, 
        and we classify 
        their finite-dimensional irreducible representations.
\end{abstract}

\vfill
\rightline{\tt math.QA/0209339}
\rightline{\LAP-938/02}
\rightline{September 2002}

\newpage
\pagestyle{plain}
\setcounter{footnote}{0}
\setcounter{tocdepth}{2}
\tableofcontents
\newpage
\section{Introduction}
\cw-algebras have been introduced in the $2d$-conformal models as a 
tool for the study of these theories. Then, these algebras and their 
finite-dimensional versions appeared to be relevant in several 
physical backgrounds. For more details on \cw-algebras, see e.g. \cite{walg}.
However, a full understanding of their algebraic 
structure (and of their geometrical interpretation) is 
lacking. The connection of some  of these finite \cw-algebras with Yangians 
appeared to be a solution at least for the algebraic structure. It 
could be surprising that Yangians \cite{Drin}, which play an important role in 
integrable systems, see e.g. \cite{ISY}, enter into the study of algebras originating  
from $2d$-conformal models. Let us however note that such a connection 
have been already remarked in WZW models \cite{Ywzw}. For more informations 
on the algebraic structure of Yangians, see e.g. \cite{Molev2001} 
and ref. therein.

The existence of an algebra 
homomorphism between Yangian based on $sl(N)$ and finite 
$\cw(sl(Np),N.sl(p))$-algebras was first proved in \cite{RS}. 
Such a connection plays a role in the study of physical models: for 
instance, in the case of the $N$-vectorial non-linear Schr\"odinger equation on the 
real line, the full symmetry is the Yangian $Y(gl(N))\equiv Y(N)$, but the space of states 
with particle number less than $p$ is a representation of the 
$\cw(gl(Np),p.sl(N))$ algebra \cite{schro}.

Later, the connection between Yangians and finite $\cw(gl(Np),N.sl(p))$-algebras 
was proven in the FRT presentation \cite{FRT} of the Yangian. It 
appears that in this framework the above $\cw$-algebras are nothing 
but truncations of the  Yangian $Y(N)$, $p$ indicating the level where 
the truncation occurs.
Thanks to this presentation, an (evaluated) $R$-matrix for these 
\cw-algebras was given, and their finite-dimensional irreducible 
representations classified \cite{wrtt}.

Then \cite{ytwist}, this connection was extended to a class of \cw-algebras, 
namely the algebras of type $\cw[so(2mp),m.sl(p)]$, 
$\cw[so((2m+1)p),m.sl(p)+so(p)]$ and $\cw[sp(2np),n.sl(p)]$, 
which where related to truncation of twisted Yangian $Y^\pm(N)$.
Note that although Yangians based on orthogonal and symplectic 
algebras exist \cite{Drin}, and admit an FRT presentation \cite{soya}, 
it is the twisted Yangians introduced by Olshanski \cite{olsh,MNO} 
which enter into the game. These later are not Hopf algebras but only 
Hopf coideals in $Y(N)$.
Nevertheless, this relation allows to give an $R$-matrix presentation 
of the $\cw$-algebras under consideration, with however the slight 
change that it is an "RSRS" relation which occurs, not an "RTT" one. 
The classification of finite-dimensional irreducible 
representations of the \cw-algebras then follows \cite{ytwist}.

\null

The aim of the present article is to extend the above correspondence 
to the case of finite \cw-superalgebras, based on Lie superalgebras 
$gl(M|N)$ and $osp(M|2n)$. As for $gl(N)$ on the one hand, and $so(m)$ 
and $sp(2n)$ on the other hand, the treatment for $gl(M|N)$ and for 
$osp(M|2n)$ will be very different.
Due to this difference, this article is divided in two main parts. 
In the first part, we show 
that $\cw(gl(Mp|Np),(M+N)gl(p))$ superalgebras are truncation of the 
super-Yangian based  on $gl(M|N)$, leading to an "RTT" presentation of 
these \cw-superalgebras.  We use this property to classify the 
finite-dimensional irreducible representations of these 
\cw-superalgebras. In the second part, we deal 
with $\cw$-superalgebras based on $osp(M|N)$ and twisted super-Yangians. 
We show that these \cw-superalgebras are truncations of twisted 
super-Yangians, leading to an "RSRS" presentation of the formers and a 
classification of their finite-dimensional irreducible 
representations.

\section{Super-Yangian}

The super-Yangian $Y(gl(M|N))=Y(M|N)$ was first defined by Nazarov \cite{Nazarov}. 
It can be obtained as the generalization of the construction for the Yangian 
$Y(M)$, 
based on the Lie algebra $gl(M)$, to the case of the Lie superalgebra $gl(M|N)$. 
Its representations have been studied by Zhang \cite{Zhang}.

\subsection{Introduction to $Y(M|N)$}

The Lie superalgebra $gl(M|N)$  
is a  $\ZZ_2$-graded vector space over $\CC$ spanned by the basis 
$\{\ce_{ab}| a,b=1,2,...,M+N\}$ . We introduce the gradation index $ [\:]$: 
\begin{displaymath}
  [a]=\left\{ \begin{array}{ll}
  0 & \mbox{if $a \leq  M$} \\
 1 & \mbox{if $M<a \leq M+N$}                      
 \end{array}\right.\mb{and}
 [\ce_{ab}]=[a]+[b]
\end{displaymath}
The bilinear graded commutator 
associated to $gl(M|N)$ is defined as follows: 
\begin{equation}
[\, , \, \} \ :\ \left\{ \begin{array}{l}
 gl(M|N)\otimes gl(M|N) \rightarrow gl(M|N) \\[1.2ex]
{[}\ce_{ab},\, \ce_{cd}\}=\delta_{cb}\ce_{ad}-(-1)^{([a]+[b])([c]+[d])}
\delta_{ad}\ce_{cb}
\end{array} \right.
\end{equation}

The super-Yangian $Y(M|N)$ is a 
$\ZZ_2$-graded Hopf algebra generated by an infinite set of elements  
$T^{ab}_{(n)}$, $a,b=1,2,...,M+N$ and  $n\in\ZZ_{>0}$. 
The $T^{ab}_{(n)}$ are
even if $[a]+[b]\equiv 0\; (mod\:2)$ and odd otherwise.

We introduce the generating function
\begin{equation}
        T(u) = \sum_{a,b=1}^{M+N}T^{ab}(u)E_{ab}\mb{and}
T^{ab}(u) = \sum_{n=0}^{\infty}T^{ab}_{(n)}u^{-n}
\end{equation}
with $T^{ab}_{(0)}=\delta^{ab}$, $u$ a spectral parameter and $E_{ab}$ 
 the matrix with 1 at position $(a,b)$ and 0 elsewhere.

The following $R$ matrix 
\begin{displaymath}
\begin{array}{cccc}
R(u)&=& \unity\otimes\unity - 
\frac{\displaystyle P}{\displaystyle u}&\;\;\;u\in\CC
\end{array}
\end{displaymath}
satisfies the graded Yang-Baxter equation. The permutation operator $P$ is
defined by:
\begin{equation}
        P_{12}=\sum_{{i,j}}(-1)^{[j]} E_{ij}\otimes E_{ji}
\end{equation}
and the tensor product is chosen graded:
\begin{equation}
(E_{ij}\otimes E_{kl})\cdot(E_{mn}\otimes E_{pq})=
(-1)^{([k]+[l])([m]+[n])}E_{ij}E_{mn}\otimes E_{kl}E_{pq}
\end{equation}
The defining relations in $Y(M|N)$ can be written as follows:
\begin{equation}
\label{RTT=TTR}
R(u-v)T_1(u)T_2(v)\;=\;T_2(v)T_1(u)R(u-v)
\end{equation}
with
\begin{equation}
 T_1(u)  = 
\sum_{a,b=1}^{M+N} T^{ab}(u)\ E_{ab}\otimes\unity \mb{and}
T_2(v)  = \sum_{a,b=1}^{M+N} T^{ab}(v)\ \unity\otimes E_{ab} 
\end{equation}

We can rewrite equation (\ref{RTT=TTR}) as follows:
\bea
[T^{ab}(u),T^{cd}(v)\}& = &\frac{(-1)^{[c]([a]+[b])+[a][b]}}{u-v}
\Big(T^{cb}(u)T^{ad}(v)-T^{cb}(v)T^{ad}(u)\Big) 
\end{eqnarray}
or equivalently:
\begin{eqnarray}
[T^{ab}_{(m)},T^{cd}_{(n)}\} & = & 
\delta^{cb}T^{ad}_{(m+n-1)}-(-1)^{([a]+[b])([c]+[d])}
\delta^{ad}T^{cb}_{(m+n-1)}+
\label{relations-Yangien}\\
&& +(-1)^{[c]([a]+[b])+[a][b]}\sum_{r=1}^{ min-1}
\left\{T^{cb}_{(r)}T^{ad}_{(m+n-1-r)}-
T^{cb}_{(m+n-1-r)}T^{ad}_{(r)}\right\}
\nonumber
\end{eqnarray}
where $min$ stands for $min(m,n)$.

The Hopf structure is given by: 
\begin{align}
&\epsilon(T^{ab}(u))=\delta^{ab} \qquad
S(T^{ab}(u))=(T^{-1}(u))^{ab}\\
&\Delta(T^{ab}(u))=\sum_{e=1}^{M+N}T^{ae}(u)\otimes T^{eb}(u) 
\label{Hopf structure}
\end{align}

The super-Yangian $Y(M|N)$ is a deformation of 
the enveloping algebra of a polynomial algebra 
(restricted to its positive modes) based on $gl(M|N)$, noted 
$\cu(gl(M|N)[x])$. 
The parameter 
$\hbar$ can be recovered after rescaling of the generators by an 
appropriate power of $\hbar$
$T^{ab}_{(n)} \rightarrow \hbar^{n-1} T^{ab}_{(n)}$.

\subsection{Finite-dimensional irreducible representations of $Y(M|N)$}

The finite-dimensional irreducible representations of $Y(M|N)$ have been 
studied in \cite{Zhang}. We recall here the main results, using a 
different basis for the positive roots (see \cite{repSY} for details).

We introduce $\NN_{M+N}=[1,M+N]\cap\ZZ_{+}$,
$\NN_{M+N}^{2}= \NN_{M+N}\times \NN_{M+N}$
and
\begin{equation}
\Phi^{+}=\Big\{(a,b)\in \NN_{M+N}^{2}, \mbox{ with either }
\left|\begin{array}{l}
1\leq a<b\leq M \\
M+1\leq a<b\leq M+2n \\
1\leq a\leq M\mbox{ and } M+n+1\leq b\leq M+2n\\
M+n+1\leq a\leq M+2n \mbox{ and }  1\leq b\leq M
\end{array}\right.
\Big\}
\end{equation}

\begin{defi}:
Let $V$ be an irreducible $Y(M|N)$-module. A nonzero element 
$v_{+}^\Lambda~\in~V$ is called {highest weight vector} if 
\begin{eqnarray}
T^{ab}_{(n)}v_{+}^\Lambda=0, & \forall  
(a,b)\in\Phi_{+} & n>0 
\nonumber \\
T^{aa}_{(n)}v_{+}^\Lambda=\lambda_a^{(n)}v_{+}^\Lambda, & a=1,..., M+N, & 
n>0,\quad
\lambda_a^{(n)}\in\CC. 
\end{eqnarray}

An irreducible module is called a {highest weight module} if it 
admits a highest weight vector.
We define  
\begin{equation}
\label{Def-Lambda}
\Lambda (u)\equiv (\lambda_1(u),\lambda_2(u),...,\lambda_{M+N}(u))
\end{equation}
with $\lambda_a(u)=1+\sum_{n>0}\lambda^n_a u^{-n}$ and call 
$\Lambda(u)$ {a highest weight} of $V$.
\end{defi}

\begin{theo}:\label{thm-Y.hw}
Every finite-dimensional irreducible $Y(M|N)$-module $V$ contains a unique 
(up to scalar multiples) highest weight vector $v_{+}^\Lambda$.

Corresponding to each $\Lambda(u)$ of the form (\ref{Def-Lambda}), 
there exists a unique irreducible highest weight $Y(M|N)$-module 
$V(\Lambda)$ with highest weight $\Lambda(u)$.
\end{theo}

\begin{theo}:
The irreducible highest weight $Y(M|N)$-module $V(\Lambda)$ is finite 
dimensional if only if its highest weight $\Lambda(u)$ satisfies
the following conditions:
\begin{equation}
\begin{array}{llll}
\displaystyle \frac{\lambda_a (u)}{\lambda_{a+1} (u)} & = & 
\displaystyle\frac{P_a (u+1)}{P_{a}(u)} & 1\leq a < N+M,\; 
a\neq M \\
&&& \\
\displaystyle\frac{\lambda_{M}(u)}{\lambda_{M+1}(u)} & = & 
\displaystyle\frac{\tilde{P}_{M}(u)}{P_{M}(u)}&
\end{array} 
\end{equation}
where, $m_a$ being the degree of $P_a$,
\begin{equation}
\begin{array}{llll}
P_{a}(u)&=&\displaystyle \prod_{i=1}^{m_{a}}(u-\gamma_{a}^{(i)})& 1\leq a 
< N+M\mbox{ and } a\neq M, \ 
\ \gamma_{a}^{(i)}\in\CC \\[1.2ex]
\tilde{P}_{M}(u)&=& \displaystyle \prod_{i=1}^{m_{M}}\left(1-
\frac{\tilde{r}^{(i)}}{u}\right)&\mb{and}
P_{M}(u)=\displaystyle  \prod_{i=1}^{m_{M}}\left(1-\frac{r^{(i)}}{u}
\right),\quad r^{(i)},\tilde{r}^{(i)}\in\CC
\end{array} 
\end{equation}
\end{theo}

Among the finite-dimensional highest weight representations, 
there is a class of particular interest: 

\begin{defi}[Evaluation representations]\hfill\\
\label{def.eval}
An evaluation representation $ev_{\pi_{\mu}}$ is a morphism from 
the super-Yangian 
$Y(M|N)$ to a highest weight irreducible representation $\pi_{\mu}$ of 
$gl(M|N)$.
The morphism is given by:
\begin{equation}
ev_{\pi_{\mu}}(T^{ab}(u))=\delta^{ab} + 
\pi_{\mu} (\ce^{ab})u^{-1}\quad \forall a,b\in\{1,...,M+N\}
\end{equation}
that is
\begin{equation}
ev_{\pi_{\mu}}(T^{ab}_{(0)})=\delta^{ab}\ ;\ 
ev_{\pi_{\mu}}(T^{ab}_{(1)})=\pi_{\mu} (\ce^{ab})
\ ;\ ev_{\pi_{\mu}}(T^{ab}_{(r)})=0\mb{for} r>1
\end{equation}
where $\ce^{ab}$ are the standard $gl(M|N)$ generators.

The highest weight $\mu(u)=(\mu_1(u),...,\mu_{M+N}(u))$ of
the representation $ev_{\pi_{\mu}}$ is given by:
\begin{equation}
\begin{array}{llll}
\mu_{a}(u)&=&1+\mu_{a} u^{-1}&\forall a\in \{1,...,M+N\}
\end{array}
\end{equation}
where $\mu = (\mu_{1},...,\mu_{M+N})$ is the highest weight of $\pi_{\mu}$.
\end{defi} 

Any finite-dimensional irreducible representation can be obtained through the
tensor products\footnote{Note however that one has sometimes to 
make a quotient to get an irreducible representation from these 
tensor products.} of such evaluation representations \cite{Zhang}:

\begin{defi}[Tensor product of evaluation representations]\hfill\\
\label{tens.eval}
Let be $\{ev_{\pi_i}\}_{i=1,...,s}$ be a set of evaluation representations. 
The tensor products of these $s$ representations 
$ev_{\vec{\pi}}=ev_{\pi_1}\otimes ...\otimes ev_{\pi_s}$ is a morphism 
from $Y(M|N)$ to the tensor product 
of $gl(M|N)$ representations $\vec{\pi}=\otimes_{i}\pi_{i}$ given by:
\begin{equation}
ev_{\vec{\pi}}(T^{ab}_{(r)})=\opluslim_{r_{1}+r_{2}+..+r_{n}=r}
\left( ev_{\pi_{1}}(T^{ai_{1}}_{(r_{1})})\otimes
ev_{\pi_{2}}(T^{i_{1}i_{2}}_{(r_{2})})\otimes\cdots\otimes
ev_{\pi_{s}}(T^{i_{s-1}b}_{(r_{s})})\right)
\end{equation}
where there is an implicit summation on the indices 
$i_{1},i_{2},\ldots,i_{s-1}=1,\ldots,M+N$.

It satisfies:
\begin{equation}
ev_{\vec{\pi}}(T^{ab}_{(r)})\neq 0 \mbox{ if and only if } r\leq s
\label{tensW}
\end{equation}
\end{defi}

\subsection{Truncated super-Yangians}
We will proceed as in \cite{wrtt}: we introduce 
$\mathcal{T}_{p}\equiv\mathcal{U}(\{T^{ij}_{(n)}, n>p\})$ and the left 
ideal $\mathcal{I}_{p}\equiv Y(M|N)\cdot\mathcal{T}_{p}$ generated by 
$\mathcal{T}_{p}$. We then define the 
coset (truncation of the super-Yangian at order $p$):
\begin{equation}
Y_{p}(M|N)\equiv Y(M|N)/\mathcal{I}_{p}
\end{equation}

\begin{prop}
The truncated super-Yangian $Y_{p}(M|N)$ is a superalgebra 
($\forall p\in\ZZ_{>0}$). 
\end{prop}
\underline{Proof}: As in \cite{wrtt} the Lie superalgebra 
structure of $Y_{p}(M|N)$ can be proved by showing that 
$\mathcal{I}_{p}$ is a two-sided ideal.
We first show that
\begin{equation}
\left[ Y(M|N)\, , \, \ct_{p}\right]\subset Y(M|N)\cdot \ct_{p}=\ci_{p}
\label{eq:ideal}
\end{equation}
The relation (\ref{relations-Yangien}) shows that 
$[T^{ij}_{m},T^{kl}_{n}]$ (for $n>p$) is the sum of two terms, 
the first being in $Y(M|N)\cdot \ct_{p}$, the second belonging to 
$\ct_{p}\cdot Y(M|N)$. Focusing on the latter, one rewrites it as
\begin{eqnarray}
\hspace{-2.4ex}&&
\sum_{r=0}^{\mu-1}\left(T^{il}_{r} T^{kj}_{m+n-1-r} +
(-1)^{[i]([k]+[j])+[k][j]}
\sum_{s=0}^{r-1}\left(T^{ij}_{s}T^{kl}_{m+n-2-s}
-T^{ij}_{m+n-2-s}T^{kl}_{s}\right)\right) 
\label{TT=TT}\\
&&\ =\  
\sum_{r=0}^{\mu-1}T^{il}_{r} T^{kj}_{m+n-1-r} +
(-1)^{[i]([k]+[j])+[k][j]}
\sum_{s=0}^{\mu-2}(\mu-s-1)\left(
T^{ij}_{s}T^{kl}_{m+n-2-s}-T^{ij}_{m+n-2-s}T^{kl}_{s}\right) 
\nonumber
\end{eqnarray}
where $\mu$ stands for min$(m,n)$. In (\ref{TT=TT}), the first sum 
belongs to $\ci_{p}$, while
 the last sum belongs to $\ct_{p}\cdot Y(M|N)$, 
with a summation which has one term less than the previous one:
we can thus proceed recursively in a finite number of steps. The final result 
is an element of $Y(M|N)\cdot \ct_{p}$.
In the same way, one can show that
\begin{equation}
\left[ Y(M|N)\, , \, \ct_{p}\right]\subset \ct_{p}\cdot Y(M|N) 
\end{equation}
so that $\ci_{p}=Y(M|N)\cdot \ct_{p}=\ct_{p}\cdot Y(M|N) $.
\finprf
Note that $\Delta$ is not a morphism of this superalgebra 
(for the structure induced by $Y(M|N)$), {\em i.e.} $Y_{p}(M|N)$ has
no natural Hopf structure.

Finally, we remark that each $Y_{p}(M|N)$ is a deformation of a truncated
polynomial algebra based on $gl(M|N)$. By truncated polynomial algebra, we mean the 
quotient  of a usual $gl(M|N)$ polynomial algebra (of generators 
$T^{ij}_{(n)}$) by the 
relations $T^{ij}_{(n)}=0$ for $n>p$. 
The construction is the same as  for the 
full super-Yangian.

\subsection{Poisson super-Yangians \label{classY}}

In the following we will deal with classical super-Yangian, where the 
commutator is replaced by a $\ZZ_{2}$-graded Poisson Bracket (PB). 
It corresponds to the usual 
classical limit of quantum groups. One sets:
\begin{eqnarray}
L(u)&=&\sum_{a,b=1}^{M+N}(-1)^{[b]}T^{ab}(u)\otimes E_{ba}\nonumber\\
R_{12}(u)&=&\unity\otimes\unity + \hbar\; r_{12}(u) 
+ o(\hbar) \mbox{ with } r_{12}(u)=
\frac{P_{12}}{u}\nonumber\\
\left[ \; , \; \right\} & = &\hbar\{\; , \;\} + o(\hbar) \nonumber
\end{eqnarray}
The relation (\ref{RTT=TTR}) is then expanded as a series in $\hbar$.
Since in a classical super-Yangian we have
$T^{ab}_{(n)}T^{cd}_{(m)}=(-1)^{([a]+[b])([c]+[d])}T^{cd}_{(m)}T^{ab}_{(n)}$, 
we obtain:
\begin{eqnarray}
\left\{T^{ab}(u),T^{cd}(v)\right\}&=&\displaystyle \frac{1}{u-v} 
(-1)^{[c]([a]+[b]) 
+[a][b]}\Big(T^{cb}(u)T^{ad}(v)-T^{cb}(v)T^{ad}(u)\Big)
\end{eqnarray}
which leads to:
\begin{eqnarray}
\{T^{ab}_{(m)},T^{cd}_{(n)}\} & = & \delta_{cb}T^{ad}_{(m+n-1)}-
(-1)^{([a]+[b])([c]+[d])}\delta_{ad}T^{cb}_{(m+n-1)}+
\label{PB-Yangian}\\
&& +(-1)^{[c]([a]+[b])+[a][b]}\sum_{r=1}^{ min(m,n)-1}
\left(T^{cb}_{(r)}T^{ad}_{(m+n-1-r)}-T^{cb}_{(m+n-1-r)}T^{ad}_{(r)}\right)
\nonumber
\end{eqnarray}

In classical super-Yangians, all the algebraic properties described above
still apply.

\section{$\cw(gl(Mp|Np),(M+N)gl(p))$ superalgebras}

For simplicity we note $\cw_{p}(M|N)\equiv\cw(gl(Mp|Np),(M+N)gl(p))$

\subsection{Definition of $\cw(\cg,\ch)$ superalgebras and Dirac brackets
\label{Dirac}}
$\cw(\cg,\ch)$ (super)algebras can be constructed as Hamiltonian 
reduction on a Lie (super)algebra $\cg$, with Poisson Brackets $\{.,.\}$.
The construction is done as follows.

We start with an $sl({2})$ embedding in $\cg$, this embedding being defined 
as the principal embedding in a regular sub(super)algebra 
$\ch\subset\cg$. We remind that the principal $sl({2})$ embedding of an 
algebra $\ch$ is given by $e_{+}=\sum_{i}e_{i}$, where $e_{+}$ is 
the positive root generator of $sl({2})$, and $e_{i}$ are the simple 
roots generators of $\ch$. If $\ch$ is a superalgebra, the principal 
$sl({2})$-embedding is defined as the  principal embedding of
its bosonic part.

Once the $sl({2})$ embedding in $\cg$ is fixed (\ie when $\ch\subset\cg$ 
is given), let $(e_{\pm},h)$ be its 
generators, one decomposes $\cg$ into $sl({2})$ representations. This 
amounts to take a $\cg$-basis of the form $J_{jm}^i$, $-j\leq m\leq j$, 
and $i$ labeling the multiplicities, with
\begin{equation}
{[e_{\pm},J_{jm}^i]}=\alpha_{jm}J_{j,m\pm1}^i,\qquad [h,J_{jm}^i]=mJ_{jm}^i
\mb{with} \alpha_{j,m}\in\CC
\end{equation}
We take $e_{\pm}=J_{1,\pm1}^0$ and $h=J_{1,0}^0$.
Then, one introduces a set of second class constrains (in Dirac 
terminology):
\begin{equation}
J_{jm}^i=\delta^{i,0}\delta_{j,1}\delta_{m,-1}\ \mb{for} m<j,\ 
\forall j,\,\forall i
\label{Dconst}
\end{equation}
This remains to set to zero all the generators but the 
$sl({2})$-highest weights ones (which are left free), and $e_{-}$ which 
is set to 1.

\textit{The $\cw(\cg,\ch)$ (super)algebra is defined as the enveloping algebra 
generated by the $sl({2})$ highest weight generators, equipped with the 
Dirac brackets associated to the constrains (\ref{Dconst}).}

We remind that the Dirac brackets can be calculated as follows.
If $\Phi=\{\phi_{\alpha}\}_{\alpha\in I}$ denotes the set of all the
above constraints, we have
\begin{equation}
\Delta_{\alpha\beta}=\{\phi_{\alpha},\phi_{\beta}\} \mbox{ is invertible: } 
\sum_{\gamma\in I}\Delta_{\alpha\gamma}\bar\Delta^{\gamma\beta}=
\delta_{\alpha}^\beta \mbox{ where } 
\bar\Delta^{\alpha\beta}\equiv(\Delta^{-1})_{\alpha\beta}
\end{equation}
The Dirac brackets are constructed as:
\begin{equation}
\{X,Y\}_{*}\sim \{X,Y\}-\sum_{\alpha,\beta\in 
I}\{X,\phi_{\alpha}\}\bar\Delta^{\alpha\beta}\{\phi_{\beta},Y\}
\ \ \forall X,Y
\end{equation}
where the symbol $\sim$ means that one has to apply the constraints on 
the right hand side {\em once the Poisson Brackets have been computed}.

\subsection{Soldering procedure}
The soldering procedure is an alternative way to compute the PB of 
$\cw(\cg,\ch)$ algebras.
We apply it to the superalgebra $gl(Mp|Np)$ with 
generators $\ce^{jm}_{ab}$, $0\leq j\leq p-1$, $-j\leq m \leq j$, 
$a,b=1,...,M+N$ 
(see appendix \ref{app-glMN}). Let $M^{jm}_{ab}$ be the $(M+N)$ 
square matrices representing the generators $\ce^{jm}_{ab}$ in the 
fundamental representation of $gl(Mp|Np)$.
Denoting $J_{jm}^{ab}$ the dual basis, we introduce
 the matrix 
\begin{equation}
\JJ\equiv\sum_{a,b=1}^{M+N}\sum_{j=0}^{p-1}\sum_{m=-j}^{j} 
J_{jm}^{ab}M^{jm}_{ab}
\end{equation}

Let us consider an infinitesimal transformation of parameters 
$\lambda^{ab}_{jm}$. For convenience we define the matrix $\lambda\equiv 
\lambda^{ab}_{jm}M^{jm}_{ab}$. 
\begin{eqnarray}
\delta_{\lambda}\JJ&\equiv&\left(\delta_{\lambda} J_{jm}^{ab}\right)\, 
M^{jm}_{ab} =[\lambda,\JJ] =\{str(\lambda\JJ),\JJ\}\\
& = & \lambda^{ef}_{rs}\,str\left(M^{rs}_{ef} M_{cd}^{tu}\right)
\{J_{tu}^{cd},J_{jm}^{ab}\}M^{jm}_{ab}
\end{eqnarray}
where summation over repeated indices is assumed. ${[.,.]}$ denotes 
the commutator of $\ZZ_{2}$-graded matrices, and $\{.,.\}$ the PB.

We ask $\JJ$ to be of the form:
\begin{equation}
\JJ|_{g.f.}=\epsilon_{-} + \sum_{a,b=1}^{M+N}\sum_{j=0}^{p-1}
W^{ab}_{j}M^{jj}_{ab}
\end{equation}
where $\epsilon_{-}$ is the $sl(2)$ negative root generator (see 
appendix \ref{app-clebsch}).
This remains to constrain the generators $J^{ab}_{jm}$ to obey
 the following second class 
constraints:
\begin{equation}
J^{ab}_{jm}=\delta_{j,1}\delta_{m+1,0}\delta^{ab}\,,\ 
\mb{for} -j\leq m<j,\ \forall\, j,\ \forall\, a,b
\label{2ndcl}
\end{equation}

We look for transformations leaving $\JJ|_{g.f.}$ 
with the same form:
\begin{equation}
\label{W-transformation}
\delta_{\lambda}\left(\JJ|_{g.f.}\right)=
[\lambda\,,\,\JJ|_{g.f.}]\;=\;\left(\delta_{\lambda}W^{ab}_{j}\right)
M^{jj}_{ab}
\end{equation}

The parameters $\lambda_{ab}^{jm}$ are constrained and only $(M+N)^2 p$
of them are free. The equation~(\ref{W-transformation}) leads to:
\begin{eqnarray}
\lambda_{j,m+1}^{ab} & = & \sum_{k,r=0}^{p-1}\sum_{l=-k}^{k}
\sum_{e=1}^{M+N}\left(\lambda_{kl}^{ae}W_{r}^{eb}\;<k,l\; ; \; r,r|jm>- 
W_{r}^{ae}\lambda_{kl}^{eb} \;<r,r\; ; \; k,l|jm> \right)\quad
\label{equations-delta}\\
\mbox{for} & & -j  \leq m \leq j-1 \nonumber\\[1.2ex]
\delta_{\lambda}W_{j}^{ab} &=& \sum_{k,r=0}^{p-1}\sum_{l=-k}^{k}
\sum_{e=1}^{M+N}\left(\lambda_{kl}^{ae}W_{r}^{eb}\;<k,l\; ; \; r,r|jj>- 
W_{r}^{ae}\lambda_{kl}^{eb} \;<r,r\; ; \; k,l|jj> \right)\ \ 
\label{deltaW}
\end{eqnarray}
where $<\cdot\;|\;\cdot>$ are real numbers defined in appendix 
\ref{app-clebsch}. All the coefficients
$\lambda_{kl}$ can be expressed in terms of the parameters
$\lambda_{k,-k}$ and 
the generators $W$, after a straightforward but tedious use of equations 
(\ref{equations-delta}). 

On the other hand we have:
\begin{equation}
\label{intermediaire}
\delta_{\lambda}W_{j}^{ab}=\lambda^{ef}_{rs}str(M_{ef}^{rs}M^{kk}_{cd})
\left\{W^{cd}_{k},W^{ab}_{j}\right\}
\end{equation}

With appendix A of \cite{wrtt} we obtain:
\begin{equation}
str\left(M^{rs}_{ef}M^{kk}_{cd}\right)=
\delta^{rk}\delta^{s,-k}\delta_{fc}\delta_{ed}(-1)^{[d]}(-1)^{k}(2k)!(k!)^{2}
\left(\begin{array}{c} p+k\\2k+1\end{array}\right)
\end{equation}
We define 
\begin{equation}
\tilde{\lambda}_{k}^{ab}\equiv(-1)^k (2k)!(k!)^{2}
\left(\begin{array}{c} p+k\\2k+1\end{array}\right)
\lambda_{k,-k}^{ab}
\end{equation}
Equation (\ref{intermediaire}) becomes:
\begin{equation}
\label{autre-deltaW}
\delta_{\lambda}W_{j}^{ab}=\sum_{k=0}^{p-1}\sum_{c,d=1}^{M+N}
(-1)^{[d]}\tilde{\lambda}_{k}^{dc}\{W_{k}^{cd},W_{j}^{ab}\}
\end{equation}

If we now compare (\ref{deltaW}) and (\ref{autre-deltaW}), 
the $\tilde{\lambda}_{k}^{ab}$ being independent from one another, 
we get $\{W_{k}^{cd},W_{j}^{ab}\}$
as a polynomial in the  $W$'s. 

\subsection{Calculation of Poisson Brackets \label{sectPB}} 
We now give two examples of PB calculations which will be needed in 
the following.

\subsubsection{Calculation of $\{W_{0}^{ab},W_{j}^{cd}\}$}

For $j=0$ equation (\ref{deltaW}) becomes:
\begin{eqnarray}
\label{deltaW0}
\delta_{\lambda}W^{ab}_{0}&=&\sum_{k=0}^{p-1}\left(\lambda_{k,-k}^{ae}
W^{eb}_{k}\:<k,-k\; ;\; k,k|0,0>-W^{ae}_{k}\lambda_{k,-k}^{eb}\:<k,k\; 
;\; k,-k|0,0>\right)\nonumber\\
&=&\frac{1}{p}\sum_{k=0}^{p-1}\left(\tilde{\lambda}_{k}^{ae}W_{k}^{eb}-
(-1)^{([a]+[e])([e]+[b])}\tilde{\lambda}_{k}^{eb}W_{k}^{ae}\right)
\end{eqnarray}

Equation (\ref{autre-deltaW}) rewrites:
\begin{eqnarray}
\label{autre-deltaW0}
\delta_{\lambda}W^{ab}_{0}&=& \sum_{k=0}^{p-1}\sum_{c,d=1}^{M+N}
(-1)^{[d]}\tilde{\lambda}_{k}^{dc}\{W_{k}^{cd},W_{0}^{ab}\}\nonumber\\
&=& \sum_{k=0}^{p-1}\sum_{c,d=1}^{M+N}(-1)^{[d]}(-1)^{1+([a]+[b])([c]+[d])}
\tilde{\lambda}_{k}^{dc}\{W_{0}^{ab},W_{k}^{cd}\}
\end{eqnarray}

Comparing the $\tilde{\lambda}_{k}^{dc}$-components of both equations, 
we obtain:
\begin{equation}
\label{pb-W0}
(-1)^{([a]+[b])([c]+[d])+[d]}\{W_{0}^{ab},W_{k}^{cd}\}=\frac{1}{p}\left(
\delta^{bc}(-1)^{([a]+[d])([d]+[c])}W_{k}^{ad}-\delta^{ad}W_{k}^{cb}\right)
\end{equation}
If we define $\hat{W}_{k}^{ab}\equiv (-1)^{[a]} W_{k}^{ab}$, $\forall\; k$,
equation (\ref{pb-W0}) becomes:
\begin{equation}
\label{PB-W0}
\{\hat{W}_{0}^{ab},\hat{W}_{k}^{cd}\}=\frac{1}{p}\left(\delta^{cb}
\hat{W}_{k}^{ad}-\delta^{ad}(-1)^{([a]+[b])([c]+[d])}\hat{W}_{k}^{cb}\right)
\end{equation}

\subsubsection{Calculation of $\{W_{1}^{ab},W_{j}^{cd}\}$}

Using the same procedure with $j=1$ we get:
\begin{eqnarray}
\delta_{\lambda}W_{1}^{ab} & = & (-1)^{1+[d]+([a]+[b])([c]+[d])}\{W_{1}^{ab},
W_{r}^{cd}\} \nonumber \\
&=&\displaystyle\frac{3}{p(p^2-1)} \sum_{k=1}^{p-1}  \frac{k(p^2-k^2)}{2k+1}
\left[\tilde{\lambda}_{k-1},W_{k}\right]_{-}^{ab}\nonumber\\
&&+\displaystyle\frac{3}{p(p^2-1)}\sum_{k=1}^{p-1}\sum_{n\geq k}^{p-1}\left[
\left[\tilde{\lambda}_{n},W_{n-k}\right]_{-},W_{k}\right]_{+}^{ab}\nonumber\\
&&+\displaystyle\frac{3}{p(p^2-1)}\sum_{k=0}^{p-1}\frac{1}{2k+1}
\sum_{n\geq k+1}^{p}\left[\left[\tilde{\lambda}_{n-1},W_{n-1-k}\right]_{+},
W_{k}\right]_{-}^{ab}\nonumber\\
&&-\displaystyle\frac{3}{p(p^2-1)}\sum_{n\geq m>k\geq 0}^{p-1}
\frac{1}{m(2k+1)}\left[\left[\left[\tilde{\lambda}_{n},W_{n-m}\right]_{-},
W_{m-1-k}\right]_{-},W_{k}\right]_{-}^{ab}\nonumber
\end{eqnarray}
where $\displaystyle\left[\tilde{\lambda}_{x},W_{y}\right]^{ab}_{\pm}
\equiv \sum_{e=1}^{M+N}\left(\tilde{\lambda}_{x}^{ae}W_{y}^{eb}\pm W_{y}^{ae}
\tilde{\lambda}_{x}^{eb}\right)$

We use $\hat{W}_{k}^{ab}\equiv (-1)^{[a]} W_{k}^{ab}$ and identify the 
$\tilde{\lambda}_{k}^{dc}$-components on both side of the equation: 

{\allowdisplaybreaks\begin{align}
\label{PB-W1}
&{\frac{p(p^2-1)}{3}\{\hat{W}_{1}^{ab},\hat{W}_{r}^{cd}\} = }
\\
&=\  \frac{(r+1)(p^2-(r+1)^2)}{2(r+1)+1}\left(\delta^{bc} 
\hat{W}_{r+1}^{ad}-(-1)^{([a]+[b])([c]+[d])}\delta^{ad}\hat{W}_{r+1}^{cb}
\right)
\nonumber\\
&\quad +\sum_{k=1}^{r}\left\{\delta^{bc} (-1)^{[e]} \hat{W}_{k}^{ae}
\hat{W}_{r-k}^{ed}-(-1)^{([a]+[b])([c]+[d])}\delta^{ad}(-1)^{[e]}
\hat{W}^{ce}_{r-k}\hat{W}^{eb}_{k}\right.
\nonumber\\
&\quad  + \left.(-1)^{[b]([c]+[d])+[c][d]}\left(\hat{W}^{ad}_{r-k}
\hat{W}^{cb}_{k}-\hat{W}_{k}^{ad} 
\hat{W}_{r-k}^{cb}\right)\right\}
\nonumber\\
&\quad +\sum_{k=0}^{r-1}\frac{r-k}{2k+1}\left\{\delta^{bc} (-1)^{[e]} 
\hat{W}_{k}^{ae}\hat{W}_{r-k}^{ed}-(-1)^{([a]+[b])([c]+[d])}
\delta^{ad}(-1)^{[e]}\hat{W}^{ce}_{r-k}\hat{W}^{eb}_{k}\right.
\nonumber\\
&\quad  + \left.(-1)^{[b]([c]+[d])+[c][d]}\left(\hat{W}^{ad}_{k}
\hat{W}^{cb}_{r-k}-\hat{W}_{r-k}^{ad}\hat{W}_{k}^{cb}\right)\right\}
\nonumber\\
&\quad -\sum_{r\geq m >k\geq 0}^{p-1}\frac{1}{m(2k+1)}\left\{
\delta^{cb}(-1)^{[e]+[f]}\hat{W}_{k}^{ae}\hat{W}_{m-k-1}^{ef}
\hat{W}_{r-m}^{fd}\right. 
\nonumber\\
&\quad  - (-1)^{([a]+[b])([c]+[d])}\delta^{ad}(-1)^{[e]+[f]}
\hat{W}_{r-m}^{ce}\hat{W}_{m-k-1}^{ef}\hat{W}_{k}^{fb}
\nonumber\\
&\quad  +(-1)^{[b]([c]+[d])+[c][d]}\left(\hat{W}_{r-m}^{ad}(-)^{[e]}\left(
\hat{W}_{m-k-1}^{ce}\hat{W}_{k}^{eb}\right)-(-1)^{[e]}\left(
\hat{W}_{k}^{ae}\hat{W}_{m-k-1}^{ed}\right)\hat{W}_{r-m}^{cb}\right.
\nonumber\\
&\quad + \left. \hat{W}_{m-k-1}^{ad}(-)^{[e]}\left(\hat{W}_{r-m}^{ce}
\hat{W}_{k}^{eb}\right)-(-1)^{[e]}\left(\hat{W}_{k}^{ae}
\hat{W}_{r-m}^{ed}\right)\hat{W}_{m-k-1}^{cb}\right.
\nonumber\\
&\quad + \left.\left. \hat{W}_{k}^{ad}(-)^{[e]}\left(\hat{W}_{r-m}^{ce}
\hat{W}_{m-k-1}^{eb}\right)-(-1)^{[e]}\left(\hat{W}_{m-k-1}^{ae}
\hat{W}_{r-m}^{ed}\right)\hat{W}_{k}^{cb}\right)\right\}
\nonumber
\end{align}}
where summation over $e$, $f$, $g=1,...,M+N$ is assumed. We remind that
\begin{equation}
\hat{W}_{j}^{ab}=(-1)^{[a]}W_{j}^{ab}
\end{equation}

\null

{\it The $\hat{W}$-basis is the one we will work on, we shall therefore omit 
the $\hat{}$ on $W$ from now on.}

\subsection{$\cw(sl(Mp|Np),(M+N)sl(p))$ superalgebras}
The $sl(2)$ principal embedding in $(M+N)gl(p)$ is indeed an 
embedding in $(M+N)sl(p)$, \ie it commutes with the $(M+N) gl(1)$ 
generators defined by $gl(p)=sl(p)\oplus gl(1)$. Moreover, considering 
these $(M+N)gl(1)$ 
subalgebras in $gl(Mp|Np)$ which commutes with $(M+N)sl(p)$, it is easy 
to see that none of its generators is affected by the constraints 
(\ref{2ndcl}), since they are highest weights. Furthermore, 
these $gl(1)$ generators, while they do not commute with 
all the constraints, weakly commute with them. By weakly, we mean 
after use of the constraints (once the PB have been computed).
Thus, their Dirac brackets  coincide 
with their original PB. This implies that these $gl(1)$ generators 
still form $gl(1)$ subalgebras in the \cw-superalgebra.

In addition, the diagonal $gl(1)$ of 
these $(M+N)gl(1)$ subalgebras, which
corresponds to the decomposition $gl(Mp|Np)=sl(Mp|Np)\oplus gl(1)$, 
 is central for the 
 original PB. Therefore, this $gl(1)$ 
generator is still central for the Dirac brackets.
In other words, one gets
\begin{eqnarray*}
&&\cw_{p}(M|N) = \cw(gl(Mp|Np),(M+N)gl(p))=\cw(gl(Mp|Np),(M+N)sl(p))\\
&&= \cw[sl(Mp|Np)\oplus gl(1),(M+N)sl(p)]=
\cu\Big(\,\cw[sl(Mp|Np),(M+N)sl(p)]\oplus gl(1)\,\Big)
\end{eqnarray*}

\section{Truncated super-Yangians and $\mathcal{W}$-superalgebras}

\subsection{$\mathcal{W}_{p}(M|N)$ as a deformation of a truncated 
polynomial algebra}
\begin{prop}
The $\mathcal{W}_{p}(M|N)$ superalgebra is a deformation of
the truncated polynomial superalgebra $gl(M|N)_{p}$.
\end{prop}
\prf
To see that the $\cw_p(M|N)$ is a deformation of a truncated 
polynomial algebra based on $gl(M|N)$, we modify the constraints to 
\begin{equation}
\JJ=\frac{1}{\hbar} \eps_{-}+\sum_{a,b=1}^{N}\sum_{j=0}^{p-1}
\sum_{0\leq m\leq j} 
J_{jm}^{ab} M^{jm}_{ab}
\end{equation}
These constraints are equivalent to the previous ones as soon as 
$\hbar\neq0$ (they correspond to a rescaling 
$J_{jm}^{ab}\rightarrow\hbar^{-m}J_{jm}^{ab}$). 
With these new constraints, the equations associated to 
the soldering procedure read:
\begin{eqnarray}
\lambda_{j,m+1}^{ab} & = & \hbar\,\sum_{k,r=0}^{p-1}\sum_{l=-k}^{k}
\sum_{e=1}^{M+N}\left(\lambda_{kl}^{ae}W_{r}^{eb}\;<k,l\; ; \; r,r|jm>- 
W_{r}^{ae}\lambda_{kl}^{eb} \;<r,r\; ; \; k,l|jm> \right) \nonumber\\
\mbox{for} & & -j  \leq m \leq j-1 \\
\delta_{\lambda}W_{j}^{ab} &=& \sum_{k,r=0}^{p-1}\sum_{l=-k}^{k}
\sum_{e=1}^{M+N}\left(\lambda_{kl}^{ae}W_{r}^{eb}\;<k,l\; ; \; r,r|jj>- 
W_{r}^{ae}\lambda_{kl}^{eb} \;<r,r\; ; \; k,l|jj> \right)\nonumber
\end{eqnarray}
This implies that the parameter 
$\lambda_{j,m}^{ab}$ behaves as $\hbar^{j+m}$.
Then, the Poisson brackets of the $W$ generators take the form:
\begin{eqnarray}
\{W_{j}^{ab},W_{\ell}^{cd}\}_{\hbar}
&=& 
\delta^{bc}W_{j+\ell}^{ad}-(-1)^{([a]+[b])([c]+[d])}\delta^{ad}
W_{j+\ell}^{cb}-\hbar P^{abcd}_{\hbar}(W)
\end{eqnarray}
where $P^{abcd}_{\hbar}(W)$, polynomial in the $W$'s,
has only positive (or null) 
powers of $\hbar$. 
This clearly shows that the $\cw_p(M|N)$ superalgebra 
is a deformation of the superalgebra generated by $W_{j}^{ab}\equiv 
J_{jj}^{ab}$ and with defining (undeformed) Poisson brackets:
\begin{eqnarray}
\{W_{j}^{ab},W_{\ell}^{cd}\}_{0} &=& \delta^{bc}W_{j+\ell}^{ad}-
(-1)^{([a]+[b])([c]+[d])}\delta^{ad}W_{j+\ell}^{cb}
\ \mbox{ if }j+\ell<p \\
&=& 0 \ \mbox{ if }j+\ell\geq p
\end{eqnarray}
One recognizes in this superalgebra a (enveloping) polynomial algebra based on 
$gl(M|N)$ quotiented by the relations $W_{j}^{ab}=0$ if $j\geq p$. In 
other words, this algebra is nothing but a truncated polynomial algebra, and 
the \cw-superalgebra is a deformation of it.
\finprf

\begin{prop}\label{propBas}
There exist two sets of generators $\{^{\pm}\bar{W}^{ab}_{j}\}_{j=0,...}$ in 
$\mathcal{W}_{p}(M|N)$ such that, 
$\forall a,b,c,d=1,...,M+N$:
\begin{eqnarray}
\label{Formule-recursion}
\forall j \geq 1\ \ \{^{\pm}\bar{W}_{1}^{ab},\:^{\pm}\bar{W}_{j}^{cd}\} 
& = & \delta^{cb}\; ^{\pm}\bar{W}_{j+1}^{ad} - 
(-1)^{([a]+[b])([c]+[d])}\delta^{ad}\; ^{\pm}\bar{W}_{j+1}^{cb}\\
& & + (-1)^{[c]([a]+[b])+[a][b]}\left(\bar{W}_{0}^{cb}\; 
{}^{\pm}\bar{W}_{j}^{ad}-\; ^{\pm}\bar{W}_{j}^{cb}\bar{W}_{0}^{ad}
\right)\nonumber\\
&&\nonumber\\
\label{Formule-recursion2}
\forall j \geq 0\ \ \{\bar{W}_{0}^{ab},^{\pm}\bar{W}_{j}^{cd}\}& = & 
\delta^{cb}\; ^{\pm}\bar{W}_{j}^{ad}-(-1)^{([a]+[b])([c]+[d])}\delta^{ad}\; 
^{\pm}\bar{W}_{j}^{cb}
\end{eqnarray}
The generators $^{\pm}\bar{W}_{j}^{ab}$ are polynomials of degree $(j+1)$ 
in the original generators $W_{j}^{ab}$ and are recursively defined by:
\begin{eqnarray}
\label{Def-W0}
\bar{W}_{0}^{ab} &\equiv& ^{+}\bar{W}_{0}^{ab}\; = \; ^{-}\bar{W}_{0}^{ab}\; 
= \; p W_{0}^{ab}\\
\label{Def-W1}
^{\pm}\bar{W}_{1}^{ab} & = & \displaystyle \pm\frac{p(p^2-1)}{6} W_{1}^{ab} 
+ \frac{p(p\pm 1)}{2} \sum_{e=1}^{M+N} (-1)^{[e]}W_{0}^{ae}W_{0}^{eb} 
\end{eqnarray}
and for $j>1$: 
\begin{equation}
\label{Def-Wj}
^{\pm}\bar{W}_{j}^{ab} = \displaystyle\sum_{n=1}^{j+1}
\sum_{|\vec{s}|=j+1-n}\; ^{\pm}\alpha_{\vec{s}}^{n,j}
\sum_{i_{1},...,i_{n-1}=1}^{M+N}(-1)^{[i_{1}]+...+[i_{n-1}]}\: 
W_{s_1}^{a i_{1}}W_{s_2}^{i_{1} i_{2}}\cdot\cdot\cdot 
W_{s_n}^{i_{n-1} b}
\end{equation}
for some numbers $^{\pm}\alpha_{\vec{s}}^{n,j}$ determined by (\ref{Formule-recursion}). 
The summation on 
$\vec{s}$ is understood as a summation on n positive (or null) integers 
$(s_1...s_n)\equiv\vec{s}$ such that $|\vec{s}|\equiv 
\sum_{i=1}^{n}s_{i}=j+1-n$.

The subsets $\{\, ^{\pm}\bar{W}_{j}^{ab}\}_{j=0,...,p-1}$ form two bases 
of $\mathcal{W}_{p}(M|N)$, the other generators $\{\, ^{\pm}
\bar{W}_{j}^{ab}\}_{j\geq p}$ being polynomials in the basis elements.
\end{prop}
\prf
As in~\cite{wrtt} the relations (\ref{Formule-recursion}) and 
(\ref{Formule-recursion2}) can be proven by recursion on $j$.
Indeed, a direct calculation shows that (\ref{Formule-recursion2}) is obeyed by (\ref{Def-Wj}) for 
any numbers $^{\pm}\alpha_{\vec{s}}^{n,j}$. Then, (\ref{Formule-recursion}) 
uniquely determine these numbers, up to the choice made in (\ref{Def-W1}).
\finprf
\begin{rmk}\rm The relations (\ref{Formule-recursion}) allow to compute 
recursively all the PB of $\cw_{p}(M|N)$ but $\{^{\pm}\bar{W}_{j}^{0} 
,^{\pm}\bar{W}_{k}^{0}\}$, where
\begin{equation}
{}^{\pm}\bar{W}_{j}^{0} =\sum_{a=1}^{M+N} {}^{\pm}\bar{W}_{j}^{aa} 
\end{equation}
In the following, we will assume that
\begin{equation}
\{{}^{\pm}\bar{W}_{j}^{0} ,{}^{\pm}\bar{W}_{k}^{0} \} =0,\ \forall\,j,k
\label{Wj0Wk0}
\end{equation}
Note that (\ref{Formule-recursion}) and 
(\ref{Formule-recursion2}) prove that (\ref{Wj0Wk0}) is valid for 
$j=0,1$ and $\forall\, k$. Let us also remark that, since 
$\cw_{p}(M|N)$ is a deformation of $gl(M|N)$ (see below), the lemma 
\ref{lemCocycle} ensures that $\{^{\pm}\bar{W}_{j}^{0} 
,^{\pm}\bar{W}_{k}^{0}\}$ is central in $\cw_{p}(M|N)$.
\end{rmk}

The first and the last coefficients that appear 
in definition (\ref{Def-Wj}) can be computed by recursion ($\forall j\geq 0$):
\begin{eqnarray}
^{\pm}\alpha^{1,j}_{j} & = & \displaystyle (\pm 1)^{j}(j!)^2
\left(\begin{array}{cc}
p+j \\ 2j+1\end{array} \right)  \label{coef1}\\
^{-}\alpha^{j,j+1}_{(0,...,0)} & = & \displaystyle 
\left(\begin{array}{cc} p \\
j+1 \end{array}\right) \label{mcoefin}\\
^{+}\alpha^{j,j+1}_{(0,...,0)}  &=& \displaystyle \left(
\begin{array}{cc} p+j \\ j+1 \end{array}\right)\label{pcoefin}
\end{eqnarray}
The non-vanishing coefficients (\ref{coef1}) show that the generators 
${}^\pm\bar W^{ab}_{j}$ for $j<p$ are indeed independent, since these 
generators write ${}^\pm\bar W^{ab}_{j}= 
{}^{\pm}\alpha^{1,j}_{j}W^{ab}_{j}+\mbox{lower}$, where lower is a 
polynomial in $W_{k}$'s with $k<j$.
\begin{coro}
The change of generators between $\{^{+}\bar{W}^{ab}_{j}\}_{j=1,...}$ and 
$\{^{-}\bar{W}^{ab}_{j}\}_{j=1,...}$ is given by: 
\begin{eqnarray}
^{\pm}\bar{W}_j^{ab} & = &\displaystyle \sum_{n=1}^{j+1}(-1)^{j+1+n}
\sum_{|\vec{s}|=j+1-n}\; \sum_{\;i_1,...,i_{n-1}=1}^{M+N}\; ^{\mp}
\bar{W}_{s_1}^{a i_1}...^{\mp}\bar{W}_{s_n}^{i_{n-1} b}\; 
(-1)^{[i_1]+...+[i_{n-1}]}
\label{eq.chgB}
\end{eqnarray}
\end{coro}
\underline{Proof}: The procedure is the same as in \cite{wrtt}:
a direct calculation shows that indeed the expression (\ref{eq.chgB}) 
satisfies (\ref{Formule-recursion}-\ref{Formule-recursion2}), and that 
(\ref{eq.chgB}) is valid for ${}^\pm\bar W^{ab}_{1}$.
\finprf

\begin{coro}
The basis $\{^{-}\bar{W}^{ab}_{j}\}_{j=1,...,p-1}$ is such that 
$^{-}\bar{W}^{ab}_{j}=0$ for $j\geq p$. In the basis 
$\{^{+}\bar{W}^{ab}_{j}\}_{j=1,...,p-1}$ all the 
$^{+}\bar{W}^{ab}_{j}$ generators $(j\geq p)$ are not vanishing.
\end{coro}
\prf 
(\ref{pcoefin}) shows that ${}^+\bar W^{ab}_{j}\neq0$ for $j\geq p$. 
Now, using (\ref{Formule-recursion}) for $j=p$, with the form 
(\ref{Def-Wj}), one gets $\alpha^{n,p}_{\vec{s}}=(-1)^nA$ with $A=0$ 
or 1. Then, (\ref{mcoefin}) shows that $A=0$ for ${}^-\bar W^{ab}_{p}$.
Finally, (\ref{Formule-recursion}) ensures that ${}^-\bar W^{ab}_{j}=0$, 
for $j>p$,
as soon as ${}^-\bar W^{ab}_{p}=0$.
\finprf

\subsection{$\mathcal{W}_{p}(M|N)$ and $Y_{p}(M|N)$}
We have shown that both $\mathcal{W}_{p}(M|N)$ and $Y_{p}(M|N)$ 
are deformations of a truncated polynomial 
superalgebra based on $gl(M|N)$.  It remains to show that these 
deformations coincide.
\begin{theo}
The $\mathcal{W}_{p}(M|N)$ superalgebra is the truncated 
super-Yangian $Y_{p}(M|N)$
\end{theo}
\prf
First, the map ${}^{-}\bar W_{j}^{ab}\to T_{j-1}^{ab}$, $\forall\, 0\leq 
j<p$, 
between basis vectors shows that $\mathcal{W}_{p}(M|N)$ and $Y_{p}(M|N)$ 
are isomorphic as vector spaces (and indeed coincide with $gl(M|N)$).
Since they are both deformations of $gl(M|N)_{p}$, we can introduce 
$\vph^W$ and $\vph^T$, the cochains associated to the deformation 
corresponding to 
$\cw_{p}(M|N)$ and $Y_{p}(M|N)$ respectively.

Now, remark that the two superalgebras have identical (in fact 
undeformed) PB on the 
couples $({}^{-}\bar W_0^{ab},{}^{-}\bar W_j^{cd})$, which proves that the 
cochains $\vph^W$ and 
$\vph^T$ coincide (in fact vanish) on these points. It is also the 
case for the couples $({}^{-}\bar W_j^{0},{}^{-}\bar W_k^{0})$, due to 
the formula (\ref{PB-Yangian}) and assumption (\ref{Wj0Wk0}).

Moreover, the property \ref{propBas} shows that the cochains 
$\vph^W$ and 
$\vph^T$ coincide on the couples 
$({}^{-}\bar W_1^{ab},{}^{-}\bar W_j^{cd})$. 
Since $\vph^W$ and 
$\vph^T$ are cocycles, this is enough (using lemma 
\ref{lemCocycle}) to prove that they are identical. 
\finprf
\subsection{Representations of $\mathcal{W}_{p}(M|N)$}
\begin{theo}\label{thm.hw}
    Any finite-dimensional irreducible representation of the 
    $\cw_{p}(M|N)$
    superalgebra is highest weight. It has a unique (up to scalar 
    multiplication) highest weight vector.
\end{theo}
\prf
An irreducible representation $\pi$ of the $\cw_{p}(M|N)$ 
superalgebra can be lifted to a representation of the whole super-Yangian by 
setting $\pi(T^{ij}_{(r)})=0$ for $r>n$. It is then obviously irreducible 
for the super-Yangian, and thus is highest weight by theorem \ref{thm-Y.hw}.
\finprf
\begin{theo}\label{thm.rep}
    {\bf Finite dimensional irreducible representations of 
    $\cw_p(M|N)$}\\    
    Any finite-dimensional irreducible representation of the 
    $\cw_{p}(M|N)$
    superalgebra is isomorphic to an evaluation representation or to the 
    subquotient of tensor product of  at most $p$ evaluation representations.
\end{theo}
\prf
By evaluation representations for $\cw_{p}(M|N)$ superalgebra, we mean the 
definitions \ref{def.eval} and \ref{tens.eval} with the change 
$T^{ab}_{r}\rightarrow\, W^{ab}_{r-1}$ (\ie the evaluation 
representations of the truncated super-Yangian). The property (\ref{tensW}) 
clearly shows that the (subquotient of) 
tensor product of $n$ evaluation representations is 
a representation of the truncated super-Yangian as soon as $n\leq p$. It 
also shows that if it is irreducible for the super-Yangian, then it is also 
irreducible for the truncated super-Yangian and that they are finite 
dimensional.

Now conversely, an irreducible representation $\pi$ of the $\cw_{p}(M|N)$ 
superalgebra can be lifted to a representation of the whole super-Yangian by 
setting $\pi(T^{ij}_{(r)})=0$ for $r>n$. It is then obviously irreducible 
for the super-Yangian, and thus is isomorphic to the (irreducible 
subquotient of) tensor product of 
evaluation representations.
\finprf
\section{Twisted super-Yangians\label{superYtw}}
Twisted super-Yangian have been introduce in \cite{repSY}. We remind 
here the main results. 

We start with the super-Yangian $Y(M|2n)$, and
 introduce the transposition $t$ on matrices:
\begin{equation}
E^t_{ab}=(-1)^{[a]([b]+1)}\theta_{a}\theta_{b}\, E_{\barb\bara}
\mb{with}\left\{\!
\begin{array}{l}
\bara=M+1-a\ \mbox{ for }1\leq a\leq M\\[1.2ex]
\bara=2M+2n+1-a\ \mbox{ for }M< a\leq M+2n
\end{array}\right.
\label{transp}
\end{equation}
where the $\theta_{a}$'s are given by
\begin{equation}
\begin{array}{l}
\theta_{a}=1\ \mbox{ for }1\leq a\leq M\\[1.2ex]
\theta_{a}=\mbox{sg}(\frac{2M+2n+1}{2}-a)\ 
\mbox{ for }M+1\leq a\leq M+2n
\end{array}\label{defTheta}
\end{equation}
Note that we have the relations
\begin{equation}
(-1)^{[a]}\, \theta_{a}\theta_{\bara}=1
\mb{and} [a]=[\bara]\ \ \forall\ a
\end{equation}
Then, we define on $Y(M|2n)$:
\begin{equation}
\tau[T(u)]=\sum_{{a,b}}\tau[T^{ab}(u)]\, E_{ab}=\sum_{{a,b}}T^{ab}(-u)\, 
E^t_{ab}
\end{equation}
which reads for the super-Yangian generators:
\begin{equation}
\tau(T^{ab}(u))=(-1)^{[a]([b]+1)}\theta_{a}\theta_{b}\ T^{\barb\bara}(-u)
\label{supertau}
\end{equation}
$\tau$ is  an algebra automorphism of $Y(M|2n)$.

One defines in $Y(M|2n)$:
\begin{eqnarray}
S(u) &=& T(u)\, \tau[T(u)]
=\sum_{a,b=1}^{M+N}S^{ab}(u)E_{ab}
=\II+\sum_{a,b=1}^{M+N}\sum_{n>0}u^{-n}S^{ab}_{(n)}E_{ab}\\
S^{ab}_{(n)} &=& \sum_{c=1}^{M+N}\sum_{p=0}^n(-1)^{p}
(-1)^{[c]([b]+1)} 
\theta_{c}\theta_{b}T^{ac}_{(n-p)}T^{\barb\barc}_{(p)}
\label{Sabn}\\
S^{ab}(u) &=& \sum_{c=1}^{M+N} (-1)^{[c]([b]+1)} 
\theta_{c}\theta_{b}T^{ac}(u)T^{\barb\barc}(-u)
\label{Sab(u)}
\end{eqnarray}
\begin{defi}
\label{twSY}
$S(u)$ defines a subalgebra of the super-Yangian, the
twisted super-Yangian $Y(M|2n)^{+}$. It obeys
 the following relation:
\begin{equation}
R_{12}(u-v)\, S_{1}(u)\, R'_{12}(u+v)\, S_{2}(v) = 
S_{2}(v)\, R'_{12}(u+v)\, S_{1}(u)\, R_{12}(u-v)
\label{rsrs}
\end{equation}
where $R(x)$ is the super-Yangian $R$-matrix,
\begin{equation}
R'(x)=\II+\frac{1}{x}Q=R^{t_{1}}(-x) \ \mbox{ with }\ Q=P^{t_{1}}\label{matRt}
\end{equation}
 and $t_{1}$ is the transposition 
(\ref{transp}) in the first auxiliary space.
\end{defi}
Introducing 
\begin{equation}
\tau(S(u))=\sum_{a,b=1}^{M+N}\,S^{ab}(-u)\,E_{ab}^t
\end{equation}
one gets 
\begin{equation}
\tau(S^{ab}(u))=(-1)^{[a]([b]+1)} \theta_{a}\theta_{b}S^{\barb\bara}(-u)
\end{equation}
Then, using the expression (\ref{Sab(u)}) and the commutation relations of 
the super-Yangian, one can show the symmetry relation:
\begin{equation}
\tau(S(u))=S(u)+\frac{\theta_{0}}{2u}(S(u)-S(-u))
\label{tauS}
\end{equation}

Note that the relation (\ref{rsrs}) is equivalent to the 
following commutator:
 \begin{eqnarray}
{[}S_{1}(u),S_{2}(v)] &=& \frac{1}{u-v}\left( \rule{0ex}{2.4ex}
P_{12}S_{1}(u)S_{2}(v) - S_{2}(v)S_{1}(u)P_{12}\right)+ \nonu
&&-\frac{1}{u+v}\left( \rule{0ex}{2.4ex}
S_{1}(u)Q_{12}S_{2}(v)-S_{2}(v)Q_{12}S_{1}(u)
\right)+ \label{comS1S2}\\
 &&+\frac{1}{u^2-v^2}\left(\rule{0ex}{2.4ex}
P_{12}S_{1}(u)Q_{12}S_{2}(v)-S_{2}(v)Q_{12}S_{1}(u)
P_{12}\right)\nonumber
\end{eqnarray}
and also to
 \begin{eqnarray}
&&{[}S^{ab}(u),S^{cd}(v)\} = \frac{(-1)^{([a]+[b])[c]}}{u-v}\,
(-1)^{[a][b]}\left( \rule{0ex}{2.4ex}
S^{cb}(u)S^{ad}(v) - S^{cb}(v)S^{ad}(u)\right)+ \nonu
&&\hspace{2.1ex}-\frac{(-1)^{([a]+[b])[c]}}{u+v}\left( \rule{0ex}{2.4ex}
(-1)^{[a][c]}\theta_{b}\theta_{\barc}S^{a\barc}(u)S^{\barb d}(v)
-(-1)^{[b][d]}\theta_{\bara}\theta_{d}S^{c\bara}(v)S^{\bard b}(u)
\right)+ \nonu
&&\hspace{2.1ex} +\frac{(-1)^{([a]+[b])[c]}}{u^2-v^2}\,
(-1)^{[a]}\theta_{a}\theta_{b}\left(\rule{0ex}{2.4ex}
S^{c\bara}(u)S^{\barb d}(v)-S^{c\bara}(v)S^{\barb d}(u) \right)
\label{comSijSkl}
\end{eqnarray}
As for $Y(M|N)$, one can show that $Y(M|2n)^+$ is a deformation of 
$\cu(osp(M|2n)[x])$.
\subsection{Finite-dimensional irreducible representations of twisted 
super-Yangians}
The finite-dimensional irreducible representations of twisted 
super-Yangians have been 
studied in \cite{repSY}. We recall here the main results. 
As for super-Yangian, they 
rely on the evaluation morphism:
\begin{prop}
 The following map defines an algebra inclusion:
 \begin{equation}
 \begin{array}{l}
Y(M|2n)^+\ \rightarrow\ \cu[osp(M|2n)]\\[1.2ex]
\displaystyle
S(u) \ \rightarrow\ \FF(u)=\II+\frac{1}{u+\half}F
\end{array}
\end{equation}
where the $osp(M|2n)$ generators $J^{ab}$ have been gathered 
 in the matrix
 \begin{equation}
F=\sum_{a,b=1}^{M+N}J^{ab} F_{ab}
\mb{with}
F_{ab}=E_{ab}-(-1)^{[a]([b]+1)}\theta_{a}\theta_{b}\ E_{\barb\bara}
 \end{equation}
\end{prop}
Using the above inclusion, one constructs from any finite-dimensional 
irreducible representation of $osp(M|2n)$, a finite-dimensional 
irreducible representation of $Y(M|2n)^+$.
\begin{theo}:\label{thm-Ytw.hw}
Every finite-dimensional irreducible $Y(M|2n)^+$-module contains a unique 
(up to scalar multiples) highest weight vector.
\end{theo}

A sufficient condition for 
the existence of irreducible finite-dimensional representations have 
been given in \cite{repSY}. It corresponds to an explicit construction 
of the representation as tensor product of $Y(M|N)$ evaluation 
representations and possibly one $osp(M|2n)$ representation (using the 
evaluation morphism).
These sufficient
 conditions were conjectured to be necessary: we will assume this 
conjecture in the following.

\subsection{Classical twisted super-Yangians}
As for super-Yangians, one can introduce a classical (Poisson 
bracket) version of twisted super-Yangians. The calculation is the 
same as in section \ref{classY}: one writes $R(u-v)=\II+\hbar r(u-v)$, $R'(u+v)=\II+\hbar r'(u+v)$,  and 
consider the terms in $\hbar$.  One gets:
\begin{eqnarray}
\{S_{1}(u),S_{2}(v)\} &=& r_{12}(u-v)S_{1}(u)S_{2}(v) - 
S_{2}(v)S_{1}(u)r_{12}(u-v)\nonu
&&+
S_{2}(v)r'_{12}(u+v)S_{1}(u)-S_{1}(u)r'_{12}(u+v)S_{2}(v)
\label{PB-Ytw}
\end{eqnarray}
In component, this reads:
\begin{eqnarray*}
\{S_{(q)1},S_{(r)2}\} &=& \sum_{s=0}^{\mu-1}\left[ \rule{0ex}{2.4ex}
P_{12}S_{(s)1}S_{(r+q-s-1)2}-
S_{(r+q-s-1)2}S_{(s)1}P_{12}\right.
\\
&+&\left.(-1)^{q+s}\left(\rule{0ex}{2.1ex}S_{(s)1}Q_{12}S_{(r+q-s-1)2}-
S_{(r+q-s-1)2}Q_{12}S_{(s)1}\right)\right]
 \nonumber
\end{eqnarray*}
with $\mu=min(q,p)$.

Let us remark that the symmetry relation (\ref{tauS}), in its 
classical form, takes the form
\begin{equation}
\tau(S(u))=S(-u)
\end{equation}
because the $T^{ab}(u)$ generators are now $\ZZ_{2}$-commuting.

\section{Folded \cw-superalgebras revisited\label{sectWfold}}
It is well-known that the $gl(M|N)$ superalgebra can be folded (using an outer
automorphism) into orthosymplectic ones (see e.g. \cite{dico}). In the same way,
folded \cw-superalgebras have been defined\footnote{Strictly speaking, it is the
folding of "affine" \cw-superalgebras that has been defined in \cite{fold}, but the
folding of finite \cw-superalgebras can be defined by the same procedure.}
 in \cite{fold}, and shown to be 
\cw-superalgebras based on orthosymplectic superalgebras.

We present here a different proof of this property, adapted to our purpose, 
and generalized to the case of the automorphisms presented in section 
\ref{superYtw}. For such a purpose, we use the Dirac bracket definition 
introduced in section \ref{Dirac}.

\subsection{Automorphism of $gl(Mp|2np)$ and $\cw_p(M|N)$}
As for the super-Yangian, one introduces an automorphism of $gl(Mp|2np)$ 
defined by
\begin{equation}
\tau(J_{jm}^{ab})=(-1)^{j+1} (-1)^{[a]([b]+1)}\theta^a\theta^b\ 
J_{jm}^{\barb,\bara}
\label{tau-glNp}
\end{equation}
where $\theta^a$ is defined in (\ref{defTheta}), and 
$\bara$ is given in (\ref{transp}). 

To prove that
$\tau$ is an automorphism of $gl(Mp|2np)$, we need the following property
of the Clebsch-Gordan coefficient, which was proved in \cite{wrtt}. 
Note that we need this property only for the algebra $gl(p)$, because 
of the decomposition $gl(Mp|2np)\sim gl(M|2n)\otimes gl(p)$ used 
here (see appendix \ref{app-glMN}).
\begin{prop}
\label{clebsch}
The Clebsch-Gordan like coefficients obey the rule:
\begin{equation}
<j,m;t,q|r,s>=(-1)^{j+t+r}\ <t,q;j,m|r,s>
\end{equation}
Note that in the above formula, it is not the $\ZZ_{2}$-grades $[j]$, 
$[t]$ or $[r]$ that are used, but really $j$, $t$ and $r$ themselves.
\end{prop}
With this property, it is a simple matter of calculation to show that 
$\tau$ defined in (\ref{tau-glNp}) is an automorphism of $gl(Mp|2np)$.
\subsection{Folding $gl(Mp|2np)$ and $\cw_p(M|2n)$}
\subsubsection{$gl(Mp|2np)$}
One considers the subalgebra Ker$(\II-\tau)$ in $gl(Mp|2np)$. It is generated
by the combinations:
\begin{equation}
K_{jm}^{ab}=J_{jm}^{ab}+\tau(J_{jm}^{ab})=J_{jm}^{ab}-(-1)^j\,
(-1)^{[a]([b]+1)}\,\theta^a\theta^b\ J_{jm}^{\barb\bara}
\end{equation}
which obey the symmetry relation
\begin{equation}
\tau(K_{jm}^{ab})=K_{jm}^{ab}\mb{\ie} 
K_{jm}^{ab}=(-1)^{j+1}\,
(-1)^{[a]([b]+1)}\,\theta^a\theta^b\ K_{jm}^{\barb\bara}
\end{equation}
Using the PB:
\[
\{J_{jm}^{ab},J_{k\ell}^{cd}\} = \sum_{r=|j-k|}^{j+k}\sum_{s=-r}^{r}
<j,m;k,\ell |r,s>\Big(\delta^{bc}\; J_{rs}^{ad} 
- (-1)^{([a]+[b])([c]+[d])}\,(-1)^{j+k+r}\,
\delta^{ad}\; J_{rs}^{cb}\Big)
\]
one can compute the commutation relations:
\begin{eqnarray*}
\{K_{jm}^{ab},K_{k\ell}^{cd}\} &=&
\sum_{r=\vert j-k\vert}^{ j+k}\sum_{s=-r}^r\ <j,m;k,\ell\vert r,s>
\,\left(
    \rule{0ex}{2.64ex}
    \delta^{bc}  K^{ad}_{rs} - (-1)^{j}\theta^a\theta^b
    (-1)^{[a]([b]+1)}  \delta^{\bara c}
 K^{\barb d}_{rs}\right.\\
&& \left.\rule{0ex}{2.64ex}
    -(-1)^{j+k+r}(-1)^{([a]+[b])([c]+[d])} \Big[\delta^{ad} K^{cb}_{rs} 
-(-1)^{j}\theta^a\theta^b(-1)^{[a]([b]+1)}\delta^{\barb d} 
     K^{c\bara}_{rs}\Big]
\right)
\end{eqnarray*}
After a rescaling of $K_{jm}^{ab}$, one
 recognizes the superalgebra $osp(Mp|2np)$.

\null

Looking at the decomposition of the fundamental of $gl(Mp|2np)$ with
respect to the principal embedding of $sl(2)$ in $(M+2n).sl(p)$ 
(see \cite{classW} and \cite{fold} for the technic used here) one
shows that the subalgebra $(M+2n).sl(p)$, generated by the $J^{aa}_{jm}$'s, 
is folded into a $(m+n).sl(p)$ (resp. $(m+n).sl(p)\oplus so(p)$) 
when $M=2m$ (resp. $M=2m+1$).

In the following, we will denote this subalgebra 
$[M.sl(p)]^\tau\oplus n.sl(p)$.

\subsubsection{$\cw_p(M|2n)$}
We are now dealing with the
enveloping algebra of $gl(Mp|2np)$, that we denote 
$\cu[gl(Mp|2np)]\equiv\cu(Mp|2np)$. 
One introduces the coset:
\[\begin{array}{ll}
\cu(Mp|2np)^+&\equiv\cu(Mp|2np)/\ck\
\mbox{ where } \ck=\cu(Mp|2np)\cdot\cl\\[1.2ex]
\mbox{ with } \cl&\mbox{ spanned by }
J_{jm}^{ab}-\tau(J_{jm}^{ab}),\ \forall\ a,b,j,m\\[1.2ex]
\cw_p(M|2n)^+&\equiv\cw_p(M|2n)/\cj\ \mbox{ where } 
\cj=\cw_{p}(M|2n)\cdot\ci \\[1.2ex]
\mbox{ with }\ci&\mbox{ spanned by }
W_{j}^{ab}-\tau(W_{j}^{ab}),\ \forall\ a,b,j
\end{array}
\]
We have the property
\begin{prop}
$\tau$ is an automorphism of $\cu(Mp|2np)$ 
provided with the Dirac brackets:
\begin{equation}
\tau\left(\{J_{jm}^{ab},J_{kl}^{cd}\}_*\right)=
\{\tau(J_{jm}^{ab}),\tau(J_{kl}^{cd})\}_*
\end{equation}
Hence, $\tau$ is also an automorphism of $\cw_{p}(M|2n)$.
\end{prop}
\prf
It is obvious that $\tau$ is an automorphism of Poisson brackets on 
$\cu(Mp|2np)$.
Moreover, due to the form of the constraints 
(\ref{2ndcl}), $\tau$ acts as
a relabeling (up to a sign) of the constraints:
\begin{equation}
\tau(\vph_\alpha)=\eps_{\alpha'}\vph_{\alpha'}\mbox{ where }
\alpha'\equiv\tau(\alpha)\mb{and} \eps_{\alpha'}=\eps_{\alpha}=\pm1
\end{equation}
which shows that $\tau(\Phi)=\Phi$. We have also
\begin{equation}
\tau(\Delta_{\alpha\beta})=\eps_{\alpha'}\eps_{\beta'}\Delta_{\alpha'\beta'}
\end{equation}
This implies that
\begin{equation}
\tau\left(\{A,\vph_{\alpha}\}\Delta^{\alpha\beta}\{\vph_{\beta},B\}\right)=
\{\tau(A),\vph_{\alpha'}\}\Delta^{\alpha'\beta'}\{\vph_{\beta'},\tau(B)\}=
\{\tau(A),\vph_{\alpha}\}\Delta^{\alpha\beta}\{\vph_{\beta},\tau(B)\}
\end{equation}
This shows that this automorphism is compatible with the set 
of constraints $\Phi$ and thus $\tau$ is an automorphism of the Dirac 
brackets. 
\finprf
\begin{coro}
   The Dirac brackets provide $\cw_p(M|2n)^+$ with an algebraic structure.
\end{coro}
\prf
We define on $\cw_p(M|2n)^+$ a bracket which is just  
the previous Dirac bracket restricted to this coset. Since
$\cw_p(M|2n)^+$ is generated by elements of the form $W+\tau(W)$, 
we have:
\begin{eqnarray*}
\{W+\tau(W),W'+\tau(W')\}_*&=&\{W,W'\}+\{\tau(W),\tau(W')\}_*+
\{\tau(W),W'\}_*+\{W,\tau(W')\}_*\\
&=&\{W,W'\}+\{\tau(W),W'\}_*+\tau\left(\rule{0ex}{1.2em}\{W,W'\}_*+
\{\tau(W),W'\}_*\right)
\end{eqnarray*}
\finprf
Indeed we have:
\begin{prop}
The  $\cw_p(M|2n)^+$ superalgebra is the 
$\cw[osp(Mp|2np), [M.sl(p)]^\tau\oplus n.sl(p)]$ superalgebra.

Above, the $[M.sl(p)]^\tau$  (resp. $n.sl(p)$) subalgebra 
is understood as subalgebra of the 
orthogonal (resp. symplectic) algebra in $osp(Mp|2np)$. 
\end{prop}
\prf
On the coset, we have 
$J_{jm}^{ab}\equiv\tau(J_{jm}^{ab})\equiv 2K_{jm}^{ab}$.
We introduce on $\cu(Mp|2np)$
\begin{equation}
2D\vph_\alpha=\vph_\alpha-\tau(\vph_\alpha)
\ ;\ 
2S\vph_\alpha=\vph_\alpha+\tau(\vph_\alpha)
\end{equation}
Since these generators satisfy $D\vph_\alpha=-\tau(D\vph_{\alpha})$ and
$S\vph_\alpha=\tau(S\vph_{\alpha})$ and are in $gl(Mp|2np)$, we have
\begin{equation}
\{S\vph_\alpha,D\vph_\beta\}\in\ci
\ \ie\ \{S\vph_\alpha,D\vph_\beta\}=0\ \mbox{ on }\cw_p(M|2n)^+
\end{equation}
Similarly we define
\begin{equation}
D\Delta_{\alpha\beta}=\{D\vph_\alpha,D\vph_\beta\}\ ;\ 
S\Delta_{\alpha\beta}=\{S\vph_\alpha,S\vph_\beta\}\ ;\ 
\end{equation}
which obey the properties:
\begin{eqnarray}
&&D\Delta_{\alpha\beta}=
\eps_{\alpha'}\eps_{\beta'}D\Delta_{\alpha'\beta'}
=-\eps_{\alpha'}D\Delta_{\alpha'\beta}
=-\eps_{\beta'}D\Delta_{\alpha\beta'} \label{Asym}\\
&&S\Delta_{\alpha\beta}=
\eps_{\alpha'}\eps_{\beta'}S\Delta_{\alpha'\beta'}=
\eps_{\alpha'}S\Delta_{\alpha'\beta}
=\eps_{\beta'}S\Delta_{\alpha\beta'}\label{sym}\\
&&\Delta_{\alpha\beta}=S\Delta_{\alpha\beta}+
D\Delta_{\alpha\beta}\ \mbox{ on }\cw_p(M|2n)^+
\end{eqnarray}
We will say that a matrix is $\tau$-antisymmetric when it satisfies
a relation like (\ref{Asym}), and $\tau$-symmetric when it obeys 
(\ref{sym}). 
$\tau$-antisymmetric matrices are orthogonal to 
$\tau$-symmetric ones:
\[
D\Delta\cdot S\Delta=0\mbox{ since }
(D\Delta\cdot S\Delta)_{\alpha\beta}=\sum_\gamma\,
D\Delta_{\alpha\gamma}S\Delta_{\gamma\beta}=\sum_{\gamma'}\,
D\Delta_{\alpha\gamma'}S\Delta_{\gamma'\beta}=-\sum_\gamma\,
D\Delta_{\alpha\gamma}S\Delta_{\gamma\beta}
\]
$S\Delta_{\alpha\beta}$ is the matrix of constraints of $osp(Mp|2np)$ 
reduced with respect to $[M.sl(p)]^\tau\oplus n.sl(p)$. 
Thus, it is invertible and the associated
Dirac brackets define the superalgebra 
$\cw(osp(Mp|2np),[M.sl(p)]^\tau\oplus n.sl(p))$. 
It remains
to show that, on $\cw_p(M|2n)^+$, the previously defined Dirac brackets 
coincide with these latter Dirac brackets. 

For that purpose, we use the form $\Delta=\Delta_0(\II+\wh{\Delta})$, given
in \cite{wrtt}, where $\Delta_0$ is an invertible $\tau$-symmetric matrix and
$\wh{\Delta}$ is nilpotent (of finite order $r$). 
Introducing the $\tau$-symmetrized and antisymmetrized
part of $\wh{\Delta}$, one deduces
\begin{equation*}
\Delta^{-1}=\Delta_0^{-1}\sum_{n=0}^r(-1)^n(S\wh{\Delta}+D\wh{\Delta})^n
=\Delta_0^{-1}\sum_{n=0}^r(-1)^n\left((S\wh{\Delta})^n+(D\wh{\Delta})^n\right)
=S\Delta^{-1}+D\Delta^{-1} 
\end{equation*}
which shows that $D\Delta$ is also invertible. 

On $\cw_p(M|2n)^+$, we have
\begin{eqnarray}
\{K^{ab}_{(m)},K^{cd}_{(n)}\}_* &=& \{K^{ab}_{(m)},K^{cd}_{(n)}\}-
\{K^{ab}_{(m)},D\vph_\alpha+S\vph_\alpha\}\Delta^{\alpha\beta}
\{D\vph_\beta+S\vph_\beta,K^{cd}_{(n)}\}\\
&=&\{K^{ab}_{(m)},K^{cd}_{(n)}\}-
\{K^{ab}_{(m)},S\vph_\alpha\}\Delta^{\alpha\beta}
\{S\vph_\beta,K^{cd}_{(n)}\}\\
&=& \{K^{ab}_{(m)},K^{cd}_{(n)}\}-
\{K^{ab}_{(m)},S\vph_\alpha\}(S\Delta^{\alpha\beta}+D\Delta^{\alpha\beta})
\{S\vph_\beta,K^{cd}_{(n)}\}
\end{eqnarray}
From the $\tau$-antisymmetry of $D\Delta^{-1}$, we get
\begin{equation}
\{.,S\vph_\alpha\}D\Delta^{\alpha\beta}\{S\vph_\beta,.\}=
\{.,S\vph_{\alpha'}\}D\Delta^{\alpha'\beta}\{S\vph_\beta,.\}=
-\{.,S\vph_\alpha\}D\Delta^{\alpha\beta}\{S\vph_\beta,.\}=0
\end{equation}
which leads to the Dirac brackets:
\begin{equation}
\{K^{ab}_{(m)},K^{cd}_{(n)}\}_* = \{K^{ab}_{(m)},K^{cd}_{(n)}\}-
\{K^{ab}_{(m)},S\vph_\alpha\}S\Delta^{\alpha\beta}
\{S\vph_\beta,K^{cd}_{(n)}\}
\end{equation}
These Dirac brackets are just the ones of the 
$\cw(osp(Mp|2np),[M.sl(p)]^\tau\oplus n.sl(p))$
superalgebra, by definition of $S\Delta$.
\finprf

\section{Folded \cw-algebras as truncated twisted Yangians}
\subsection{Classical case}
We start with the $\cw_p(M|2n)$ superalgebra in the Yangian basis. The Poisson
brackets are
\begin{equation}
\{T_{(m)1},T_{(n)2}\} = \sum_{r=0}^{min(m,n)-1} (
P_{12}T_{(r)1}T_{(m+n-r)2}-T_{(r)2}T_{(m+n-r)1}P_{12})
\end{equation}
with the convention $T_{(m)}=0$ for $m>p$. The action of the automorphism
$\tau$, both for twisted super-Yangian and folded $\cw_p(M|2n)$ superalgebra, reads
\begin{equation}
\tau(T_{(m)})=(-1)^{m}T^t_{(m)}
\end{equation}
However, from the twisted super-Yangian point of view, one selects
 the generators 
 \[
 S_{(m)}=\sum_{r+s=m}(-1)^{s}T_{(r)}T^t_{(s)}
 \]
while in the folded \cw-superalgebra case, one constrains the generators to
$T_{(m)}=(-1)^{m}T^t_{(m)}$. Although the procedures are different 
(and indeed lead to different generators), we have:
\begin{theo}
As an algebra, the \cw-superalgebra 
$\cw(osp(Mp|2np),[M.sl(p)]^{\tau}\oplus n.sl(p))$ is
isomorphic to the truncation (at level $p$) of the 
(classical) twisted super-Yangian $Y(M|2n)^+$. 

More precisely, we have the correspondences:
\begin{equation}
\begin{array}{ccl}
Y_{p}(2m+1|2n)^+ &\longleftrightarrow & 
\cw[osp(2mp+p|2np),(m+n).sl(p)\oplus so(p)] \\[1.2ex]
Y_{p}(2m|2n)^+ &\longleftrightarrow & \cw[osp(2mp|2np),(m+n).sl(p)]
\end{array}
\end{equation}
\end{theo}
\prf
We prove this theorem by showing that the Dirac brackets of the folded 
\cw-superalgebra coincide with the Poisson brackets (\ref{PB-Ytw})
with the truncation $S_{(m)}=0$ for $m>p$. 

We start with the $\cw_p(M|2n)$ superalgebra in the truncated super-Yangian basis:
\begin{eqnarray}
&&\{T_{(q)1},T_{(r)2}\}= \sum_{s=0}^{\mu-1}\left( P_{12}T_{(s)1}T_{(r+q-s-1)2}-
T_{(r+q-s-1)2}T_{(s)1}P_{12}\right)\nonu 
&&\mbox{ with }\ T_{(s)}=0\mbox{ for }s>p
\mbox{ and }\mu=\mbox{min}(q,r,p)
\end{eqnarray}
 In this basis, we define 
\begin{equation}
2\vph_{(s)}= T_{(s)}-(-1)^s\, T^t_{(s)} \ \mbox{ and }\ 
2K_{(s)}=T_{(s)}+(-1)^s\, T^t_{(s)}
\end{equation}
The folding (of the \cw-superalgebra) corresponds to
\begin{equation}
\vph_{(s)}=0\ \ie K_{(s)}=(-1)^sK^t_{(s)}\label{eq:fold}
\end{equation}
It is a simple matter of calculation to get
\begin{eqnarray}
2\{K_{(q)1},K_{(r)2}\} &=& \sum_{s=0}^{\mu-1}\left[ \rule{0ex}{2.4ex}
P_{12}K_{(s)1}K_{(r+q-s-1)2}-
K_{(r+q-s-1)2}K_{(s)1}P_{12}+\right.\nonu
&&\left.+(-1)^{q+s}\left(\rule{0ex}{2.1ex}K_{(s)1}Q_{12}K_{(r+q-s-1)2}-
K_{(r+q-s-1)2}Q_{12}K_{(s)1}\right)\right]
\end{eqnarray}
that is to say
\begin{equation*}
2\{K_{1}(u),K_{2}(v)\} = [r_{12}(u-v), K_{1}(u)K_{2}(v)] +
K_{2}(v)r'_{12}(u+v)K_{1}(u)-K_{1}(u)r'_{12}(u+v)K_{2}(v)
\end{equation*}
These PB are equivalent to the relation (\ref{PB-Ytw}) for 
$S(u)\equiv K(\frac{u}{2})$. 
The constraint (\ref{eq:fold}) is then rewritten as
$\tau(S(-u))=S(u)$. Thus, the folded \cw-superalgebra and the truncated 
twisted super-Yangian are defined by the same relations.
\finprf

\subsection{Quantization and representations of \cw-superalgebras}
Now that folded \cw-superalgebras have proved to be truncation of twisted 
super-Yangians at classical level, there quantization is very simple. 
It can be identified with the truncated twisted super-Yangian at 
quantum level:
\begin{eqnarray}
&&R_{12}(u-v)\, S_{1}(u)\, R'_{12}(u+v)\, S_{2}(v) = 
S_{2}(v)\, R'_{12}(u+v)\, S_{1}(u)\, R_{12}(u-v)\\[2.1ex]
&&\mbox{with }\left\{\begin{array}{l}\displaystyle
 R_{12}(x)=\II-\frac{1}{x}P_{12}\ ;\ R'_{12}(x)=\II-\frac{1}{x}Q_{12}\\
\displaystyle S(u)=\sum_{m=0}^{p}\, u^{-m}\, S_{(m)}\ ;\ S_{(0)}=\II
\end{array}\right.
\end{eqnarray}
Using the representations classification of twisted super-Yangians 
given in \cite{repSY}, one can then deduce the classification of 
irreducible finite-dimensional representations for truncated twisted 
super-Yangians in the same way it has been done in \cite{ytwist} for 
ordinary twisted Yangians. For conciseness, we will just sketch the 
results.
In particular, one gets the following 
theorems
\begin{theo}
    Any finite-dimensional irreducible representation of the 
    $\cw_{p}(M|2n)^+$
    superalgebra is highest weight.
\end{theo}
\prf
Same proof as for theorem \ref{thm.hw}.\finprf
\begin{theo}
    Any finite-dimensional irreducible representation of the 
    $\cw_{p}(M|2n)^+$
    superalgebra is isomorphic to an evaluation representation, or to the 
    (irreducible subquotient of) tensor product of at most $[p/2]$ 
    evaluation representations of $Y(M|2n)$, and possibly an $osp(M|2n)$ 
    representation.
\end{theo}
\prf
Same proof as for twisted Yangians, see \cite{ytwist}, using the 
results given in \cite{repSY} for $Y(M|2n)^+$. Indeed, as for the 
$gl(M|N)$ case, one needs to have $S_{(r)}=0$ for $r\geq p$ to get  a 
representation of the \cw-superalgebra. This constrains the number of 
evaluation representations allowed to be tensorised (to get a 
representation). The difference 
with the $Y(M|N)$ case lies in the quadratic form 
$S(u)=T(u)\tau(T(-u))$, which lowers the maximum number of terms in 
the tensor product. The occurrence of an $osp(M|2n)$ representation is 
due to the classification given in \cite{repSY}.
\finprf 
Reasoning as in \cite{ytwist}, one can also get a 
condition on the weights of the representation. We omit it here, due 
to the lack of place.
\begin{rmk}\rm As for \cw-algebras based on $so(M)$ and $sp(2n)$, 
see \cite{ytwist} for more details, one could think that 
\cw-superalgebras based on $osp(M|2n)$ are related to super-Yangians based 
on $osp(M|2n)$ instead of twisted super-Yangians. However, a simple 
counting (using the method given in \cite{classW})
of the generators shows that it is the twisted super-Yangians 
that have to be considered.
\end{rmk}

\appendix

\section{General settings on $gl(Mp|Np)$\label{app-glMN}}
\subsection{Clebsch-Gordan like coefficients\label{app-clebsch}}
We start with the $gl(Mp|Np)$ superalgebra in its fundamental representation, 
and consider
the $sl(2)$ principal embedding in $(M+N)gl(p)\equiv 
\underbrace{gl(p)\oplus ...\oplus gl(p)}_{M+N}$.

In the fundamental representation, one can view $gl(Mp|Np)$ as 
$gl(p)\,\otimes\, gl(M|N)$, so that the generators of this $sl(2)$ 
can be written as $\epsilon_{\pm,0}\equiv 
e_{\pm,0}\otimes\unity_{M+N}$. The $e_{\pm,0}$ are 
the generators  of the $sl(2)$ algebra principal in $gl(p)$
and verify  $[e_{0},e_{\pm}]=\pm e_{\pm}$ and $[e_{+},e_{-}]=e_{0}$.
The generator $\unity_{M+N}$ is the identity  generator in $gl(M|N)$. 

Under the adjoint action of this $sl(2)$, 
$gl(p)\,\otimes\, gl(M|N)$ can be decomposed in $sl(2)$ multiplets:
$M_{jm}^{ab}\equiv M_{jm}\otimes E^{ab}$
with $a,b=1,...,M+N\: ; \: -j\leq m\leq j\: ; \:  0\leq j\leq p-1$.

The $M^{jm}$ are $p\times p$ matrices resulting from the decomposition of 
$gl(p)$ in $sl(2)$ multiplets. 
Properties of the $M_{jm}$ are gathered in appendix A
of \cite{wrtt}.

The $E_{ab}$ are $(M+N)\times(M+N)$ matrices with $1$ at position $(a,b)$. 
They are the graded part of $ M_{jm}^{ab}$ which  is even if
$a+b\equiv 0 \; (mod\:2)$ and odd otherwise. 

Following appendix A of \cite{wrtt} we have
\begin{eqnarray}
[\epsilon_{+},M_{jm}^{ab}] & = & \displaystyle \frac{j(j+1)-m(m+1)}{2} 
M_{j,m+1}^{ab} \\
\;[\epsilon_{-},M_{jm}^{ab}] & = & M_{j,m-1}^{ab} \\
\;[\epsilon_{0},M_{jm}^{ab}] & = & m\; M_{jm}^{ab} 
\end{eqnarray}

The product law (in the fundamental representation) reads:
\begin{equation}
M_{jm}^{ab}\cdot M_{ln}^{cd}=\delta^{bc} \sum_{r=|j-l|}^{j+l}
\sum_{s=-r}^{r}<j,m;l,n|r,s>\; M_{rs}^{ad} 
\end{equation}
which leads to the following commutation relations (valid 
in the abstract algebra):
\begin{eqnarray*}
[M_{jm}^{ab},M_{ln}^{cd}] &=& \sum_{r=|j-l|}^{j+l}\sum_{s=-r}^{r}
\Big(\delta^{bc}<j,m;l,n|r,s>\; M_{rs}^{ad} 
\\
&& - (-1)^{([a]+[b])([c]+[d])}
\delta^{ad}<l,n;j,m|r,s>\; M_{rs}^{cb}\Big)
\end{eqnarray*}

The scalar product is:
\begin{equation}
\eta_{j,m;l,n}^{ab,cd}=str(M_{jm}^{ab}\cdot M_{ln}^{cd})=(-1)^{[a]}\delta^{ad}
\delta^{cb}(-1)^m\delta_{j,\ell}\delta_{m+n,0}\eta_{j}
\end{equation}
for some non-vanishing coefficient $\eta_{j}$, 
 given in \cite{wrtt}.

The "Clebsch-Gordan like" coefficients are then given by
\begin{equation}
<jm,k\ell |r,s>=
\frac{(-1)^s}{\eta_{r}}tr\Big(M_{jm}M_{k\ell}M_{r,-s}\Big)
\label{defClebsch}
\end{equation}
We remind that in (\ref{defClebsch}), it is the usual trace operator which is 
involved, since we are in the $gl(p)$ Lie algebra.

\subsection{Structure constants}

We consider a Lie superalgebra $\mathcal{G}$  
in its fundamental representation, with homogeneous generators $t_{a}$. 
As usual, we can define a gradation index $ [\:]$ 
such that:
\begin{displaymath}
  [a]=\left\{ \begin{array}{ll}
 0 & \mbox{if $t_a$ bosonic} \\
1 & \mbox{if $t_a$ fermionic} 
\end{array}\right.
\end{displaymath}

The commutation relations are $[t_a,t_b\}={f_{ab}}^{c}t_c$
(summation over repeated indices).
The structure constants have following property:
${f_{ab}}^{c}\neq 0 \Rightarrow [a]+[b]+[c]=0$. They obey the graded 
Jacobi identity:
\begin{equation}
{f_{ab}}^{d}{f_{dc}}^{e}={f_{bc}}^{d}{f_{ad}}^{e}+(-1)^{[b][c]}{f_{ac}}^{d}{f_{db}}^{e}
\end{equation}

Note that the adjoint representation for superalgebras takes the form
\begin{equation}
\mbox{ad}(t_{a})_{b}^c=-(-1)^{[a][b]}{f_{ab}}^c
\end{equation}

The invariant metric $g_{ab}$ is proportional to  $str_{F}(t_at_b)$, where 
the supertrace is taken in the fundamental representation. Note that 
the Killing form, which is the supertrace in the adjoint 
representation, can be degenerate (in fact null) for some
superalgebras, e.g. $gl(M|M)$ \cite{dico}. 
The invariant metric has following properties:
\begin{equation}
g_{ab}=(-1)^{[a][b]}g_{ba} \mb{and}
g_{ab}=0 \mbox{ if }[a]\neq [b]
\end{equation}

We introduce its inverse $g^{ab}$ and use it 
 to rise and lower the indices. For instance
 $t^a\equiv g^{ab} t_{b}$ and
 ${f^{ab}}_c\equiv 
g^{a\alpha}g^{b\beta}{f_{\alpha\beta}}^{\gamma}g_{\gamma c}$: 
we therefore have $[t^a,t^b]={f^{ab}}_c t^c$.

Defining the tensor $f_{abc}={f_{ab}}^\gamma\, g_{\gamma c}$, one can 
take it totally (graded) antisymmetric:
\begin{equation}
f_{abc}=-(-1)^{[a][b]}\, f_{bac}=-(-1)^{[b][c]}\, f_{acb}
\end{equation}

\section{Deformations and cohomology}

Let us consider a Lie superalgebra $\mathcal{A}$ with homogeneous 
generators $u_{\alpha}$ and Lie bracket:
\begin{equation}
\{u_{\alpha},u_{\beta}\} = {f_{\alpha\beta}}^{\gamma} u_{\gamma}
\label{PB-ca}
\end{equation}
The gradation index $[\;]$ is such that $[\alpha]=0$ if $u_{\alpha}$ is bosonic
and $[\alpha]=1$ if $u_{\alpha}$ is fermionic. 

We aim to construct a deformation of the Lie bracket (\ref{PB-ca}), 
following e.g. \cite{Gerstenhaber}.

For such a purpose, we introduce $n$-cochains ($n\in\ZZ_{>0}$), 
{\em i.e.} linear maps $\chi^{(n)}$ from $\mathcal{A}^{n}$ to $\mathcal{A}$ with 
following property: 
\begin{equation}
\chi^{(n)}(u_{\alpha_{1}},...,u_{\alpha_{i}},u_{\alpha_{i+1}},...,
u_{\alpha_{n}}) = (-1)^{1+[\alpha_{i}][\alpha_{i+1}]} 
\chi^{(n)}(u_{\alpha_{1}},...,u_{\alpha_{i+1}},u_{\alpha_{i}},...
,u_{\alpha_{n}}).
\end{equation}
The Chevalley derivation $\delta$ maps $n$-cochains to $(n+1)$-cochains:
\begin{eqnarray}
(\delta\chi^{(n)})(u_{\alpha_{0}},u_{\alpha_{1}},...,u_{\alpha_{n}}) & = & 
\sum_{i=0}^{n} (-1)^{i + \epsilon_{i}}\{u_{\alpha_{i}},
\chi^{(n)}(u_{\alpha_{0}},...,\hat{u}_{\alpha_{i}},...,u_{\alpha_{n}})\} \\
&&+ \sum_{0\leq i<j\leq n} (-1)^{i+j+\epsilon_{ij}}
\chi^{(n)}(\{u_{\alpha_{i}},u_{\alpha_{j}}\},u_{\alpha_{0}},...,
\hat{u}_{\alpha_{i}},...,\hat{u}_{\alpha_{j}},...,u_{\alpha_{n}})\nonumber
\end{eqnarray}
where $\epsilon_{i}=[\alpha_{i}](\sum_{k<i} [\alpha_{k}])$ and 
$\epsilon_{ij}=\epsilon_{i}+\epsilon_{j}+[\alpha_{i}][\alpha_{j}]$.

It obeys $\delta^{2}=0$, so that one can define $n$-cocycles, which are 
closed $n$-cochains ($\delta \chi^{(n)}=0$), and coboundaries, which are 
exact $n$-cochains ($\chi^{(n)}=\delta \chi^{(n-1)}$). As usual, one 
considers closed cochains modulo exacts ones to study the cohomology 
associated to $\delta$.

Here, we will be mainly concerned with the action of the Chevalley derivation  on 2-cochains:
\begin{eqnarray*}
(\delta\chi)(u,v,w) &=&
\{u,\chi(v,w)\}-(-1)^{[u][v]}\,\{v,\chi(u,w)\}+
(-1)^{[w]([u]+[v])}\,\{w,\chi(u,v)\}\\
&& -\chi(\{u,v\},w)+
(-1)^{[v][w]}\,\chi(\{u,w\},v)-(-1)^{[u]([v]+[w])}\,\chi(\{v,w\},u)
\end{eqnarray*}

\null

We now consider a deformation of the enveloping algebra 
$\mathcal{U}(\mathcal{A})$:
\begin{equation}
\{u_{\alpha},u_{\beta}\}_{\hbar} = {f_{\alpha\beta}}^{\gamma} u_{\gamma} + 
\hbar\varphi_{\hbar}(u_{\alpha},u_{\beta})
\end{equation}
where $\varphi_{\hbar}$ is a 2-cochain which may depend on positive 
powers of $\hbar$.  
Asking the bracket $\{\cdot, \cdot\}_{\hbar}$ to obey the
graded Jacoby identity is equivalent to say that $\vph_{\hbar}$ is a 
2-cocycle:
\begin{equation}
\delta\varphi_\hbar (u_{\alpha},u_{\beta},u_{\gamma}) = 0 
\end{equation}

We now prove a result that is used in the present article.
\begin{lem} \label{lemCocycle}
Let $gl(M|N)_{p}$ be the polynomial algebra based on 
$gl(M|N)$, truncated at order $p$, and $u_j^{ab}$ ($j<p$ and 
$a,b=1,\ldots,M+N$) the 
corresponding generators. Let $\vph$ be a 2-cocycle with values in 
$\cu(gl(M|N)_{p})$. We introduce $u_{j}^{0}=\sum_{a=1}^{M+N} u_{j}^{aa}$.

 If $\vph(u_0^{ab},u_j^{cd})$ and $\vph(u_1^{ab},u_j^{cd})$, 
$\forall a,b,c,d=1,\ldots,M+N$ and $\forall j=0,\ldots,p-1$ are 
known, then $\vph$ 
is completely determined up to $\vph(u_{j}^{0},u_{k}^{0})$, $j,k>1$,
which is 
central in $\cu(gl(M|N)_{p})$. 
\end{lem}
\prf 
We write the cocycle condition for a triplet 
$(u_j^{ab},u_k^{cd},u_\ell^{eg})$:
\begin{eqnarray}
&&\vph(\{u_{j}^{ab},u_{k}^{cd}\},u_{\ell}^{eg})+
(-1)^{([c]+[d])([e]+[g])}
\vph(\{u_{j}^{ab},u_{\ell}^{eg}\},u_{k}^{cd})\nonumber\\
&&-(-1)^{([a]+[b])([c]+[d]+[e]+[g])} 
\vph(\{u_{k}^{cd},u_{\ell}^{eg}\},u_{j}^{ab})
\ =\ \{u_{j}^{ab}, \vph(u_{k}^{cd},u_{\ell}^{eg}) \}
\label{eq-cocy}\\
&&-(-1)^{([a]+[b])([c]+[d])}\{u_{k}^{cd}, 
\vph(u_{j}^{ab},u_{\ell}^{eg}) 
\}-(-1)^{([a]+[b])([c]+[d]+[e]+[g])}\{u_{\ell}^{eg}, 
\vph(u_{j}^{ab},u_{k}^{cd}) \}
\nonumber
\end{eqnarray}
We write the commutation relations of $gl(M|N)_{p}$ as:
\begin{equation}
\{u_{j}^{ab},u_{k}^{cd}\}= \delta^{bc}\, u_{j+k}^{ad}
-(-1)^{([a]+[b])([c]+[d])}\,\delta^{ad}\, u_{j+k}^{cb}\mb{with} 
u_{n}^{ab}=0,\ n>p,\ \forall a,b
\label{PBglNp}
\end{equation}
Taking as special case $e=g=a\neq b$ and $\ell=1$, one gets from 
(\ref{eq-cocy}):
\begin{eqnarray}
\vph(u_{j+1}^{ab},u_{k}^{cd}) &=&\{u_{j}^{ab}, \vph(u_{k}^{cd},u_{1}^{aa}) \}
-\vph(\{u_{j}^{ab},u_{k}^{cd}\},u_{1}^{aa}) \label{eq-cocy-1}
\\
&+&\!(-1)^{([c]+[d])([a]+[b])}\,\Big(( \delta^{da}-\delta^{ca})
\vph(u_{k+1}^{cd},u_{j}^{ab})
-\{u_{1}^{aa},\vph(u_{j}^{ab},u_{k}^{cd})\}\nonu
&&\hspace{8.2em}
-\{u_{k}^{cd}, \vph(u_{j}^{ab},u_{1}^{aa}) \}\Big)
\nonumber
\end{eqnarray}
Taking as special case $j=1$, $k=2$, this last equation shows that one can 
compute $\vph(u_{2}^{ab},u_{2}^{cd})$, for $a\neq b$, as soon as one 
knows $\vph(u_{1}^{cd},u_{j}^{eg})$, $\forall c,d,e,g$, $\forall j$. Then, 
in the same way, $j=1$ allows to compute $\vph(u_{2}^{ab},u_{k+1}^{cd})$
as soon as one knows $\vph(u_{2}^{ab},u_{k}^{cd})$. 

More generally,  if one 
supposes by induction that $\vph(u_{j'}^{ab},u_{k}^{cd})$, 
$\forall j'\leq j$, $\forall k$,  and  $\forall\ c,d,\,a\neq b$,
are known, (\ref{eq-cocy-1}) shows that one can compute 
$\vph(u_{j+1}^{ab},u_{k}^{cd})$, for $a\neq b$ and $\forall k$. 

Thus, by induction, we have shown that one can compute 
$\vph(u_{j}^{ab},u_{k}^{cd})$, for $a\neq b$, $\forall j,k,\ c,d$,
 from the knowledge of  
$\vph(u_{1}^{cd},u_{k}^{eg})$. 

It remains to compute $\vph(u_{j}^{aa},u_{k}^{bb})$. For such a 
purpose, we start again with (\ref{eq-cocy}) with  now $a=d\neq b=c$ 
and $e=g$:
\begin{eqnarray}
\vph\big(u_{j+k}^{aa}-(-1)^{[a]+[b]}u_{j+k}^{bb}\,,u_{\ell}^{ee}\big) &=&
( \delta^{ae}-\delta^{be})\Big(
\vph(u_{j+\ell}^{ab},u_{k}^{ba})+(-1)^{[a]+[b]}
\vph(u_{k+\ell}^{ba},u_{j}^{ab}) \Big)
\nonumber\\
&&
+\{u_{j}^{ab},\vph(u_{k}^{ba},u_{\ell}^{ee})\}
-(-1)^{[a]+[b]}\{u_{k}^{ba},\vph(u_{j}^{ab},u_{\ell}^{ee}) \}\nonu
&&-(-1)^{[a]+[b]} \{u_{\ell}^{ee}, \vph(u_{j}^{ab},u_{k}^{ba}) \}
\label{eq-cocy-2}
\end{eqnarray}
All the terms in the r.h.s. of the above equation are known, so that 
one can compute\footnote{One should take $a\neq b$, but for $a=b$ one 
obviously gets 0.}
$\vph\big(\,(-1)^{[a]}u_{j}^{aa}-(-1)^{[b]}u_{j}^{bb}\,,u_{k}^{ee}\,\big)$, $\forall 
a,b,e$, $\forall j,k$.

Thus, only $\vph(u_{j}^{0},u_{k}^{0})$, where 
$u_{j}^0=\sum_{a=1}^{M+N} u_{j}^{aa}$, remains to be computed.

Once again, from (\ref{eq-cocy}), taking $a=b$ and $c=d$, and then 
 summing over $a$ and $d$, one gets
\begin{equation}
\{u_{\ell}^{eg},\vph(u_{j}^{0},u_{k}^{0})\}=0
\end{equation}
which shows  that 
$\vph(u_{j}^{0},u_{k}^{0})$ is central in 
$\cu(gl(M|N)_{p})$.

Thus, apart from the values $\vph(u_{0}^{ab},u_{k}^{cd})$ and the 
just mentioned
central terms, we are able 
to compute all the expressions $\vph(u_{j}^{ab},u_{k}^{cd})$. This ends 
the proof.
\finprf


\end{document}